\documentclass[leqno,12pt]{article}
\usepackage{latexsym}
\usepackage{amsmath}
\usepackage{amssymb}
\usepackage{bbm}
\usepackage{color}
\usepackage[authoryear]{natbib}
\usepackage{graphicx}
\usepackage{epsfig}
\newtheorem{thm}{Theorem}[section]
\newtheorem{lem}[thm]{Lemma}
\newtheorem{prop}[thm]{Proposition}
\newtheorem{cor}[thm]{Corollary}

\newtheorem{rem}[thm]{Remark}

\makeatletter\@addtoreset{equation}{section}\makeatother
\title{Sum rules via large deviations}
\author{{\small Fabrice Gamboa}\\
{\small Universit\'e Paul Sabatier}\\
{\small Institut de Math\'ematiques de Toulouse}\\
{\small 118 route de Narbonne}\\
{\small 31062 Toulouse Cedex 9, France}\\
{\small e-mail: gamboa@math.univ-toulouse.fr}\\ 
\and{\small Jan Nagel}\\
{\small Technische Universit\"at M\"unchen}\\
{\small Fakult\"at f\"ur Mathematik}\\
{\small Boltzmannstr. 3}\\
{\small 85748 Garching, Germany}\\
{\small e-mail: jan.nagel@tum.de}\\ 
\and{\small Alain Rouault}\\
{\small Universit\'e Versailles-Saint-Quentin}\\
{\small LMV UMR 8100}\\
{\small 45 Avenue des Etats-Unis}\\
{\small 78035-Versailles Cedex France}\\
{\small e-mail: alain.rouault@uvsq.fr}}
\begin{document}
\maketitle


\newcommand {\eref}[1]{(\ref{#1})}
\newcommand{\ea}{\end{array}}
\newcommand{\beqohne}{\begin{eqnarray*}}
\newcommand{\eeqohne}{\end{eqnarray*}}
\newcommand{\beohne}{\begin{equation*}}
\newcommand{\eeohne}{\end{equation*}}
\newcommand{\R}{\mathbb{R}}
\newcommand{\N}{\mathbb{N}}
\newcommand{\T}{\mathbb{T}}
\newcommand{\E}{\mathbb{E}}
\def\proof{\noindent{\bf Proof:}\hskip10pt}
\def\QED{\hfill\vrule height 1.5ex width 1.4ex depth -.1ex \vskip20pt}
\def \sur#1#2{\mathrel{\mathop{\kern 0pt#1}\limits^{#2}}}
\def \el{\sur{=}{(d)}}
\newcommand{\indi}{\mathbbm{1}} 
\newcommand{\la}{\lambda}
\newcommand{\lai}{\la_i}
\newcommand{\laj}{\la_j}
\newcommand{\mun}{\mu^{(n)}}
\newcommand{\muun}{\mu^{(n)}_{\u}}
\newcommand{\mutn}{\tilde{\mu}^{(n)}}

\newcommand{\ga}{\gamma}
\newcommand{\gai}{\ga_i}
\newcommand{\gaj}{\ga_j}
\def \u{{\tt u}}
\def \w{{\tt w}}
\def \ben{\begin{eqnarray}}
\def \een{\end{eqnarray}}
\newcommand{\tr}{\mathrm{tr}}
\newcommand{\Dir}{\mathrm{Dir}}
\newcommand{\Q}{\mathbb{Q}}
\newcommand{\Pnv}{\mathbb P^n_V}
\newcommand{\munI}{\mu^{(n)}_I}
\newcommand{\muunI}{\mu^{(n)}_{\u,I}}
\newcommand{\mutnI}{\tilde{\mu}^{(n)}_I}
\newcommand{\munIj}{\mu^{(n)}_{I(j)}}
\newcommand{\muunIj}{\mu^{(n)}_{\u,{I(j)}}}
\newcommand{\bmuunIj}{\bar{\mu}^{(n)}_{\u,{I(j)}}}
\newcommand{\mutnIj}{\tilde{\mu}^{(n)}_{I(j)}}
\newcommand{\sn}{^{(n)}}
\newcommand{\Ir}{\mathcal{J}}
\newcommand{\Fr}{\mathcal{F}}
\newcommand{\Sr}{\mathcal{S}}
\newcommand{\ap}{\alpha^+}
\newcommand{\am}{\alpha^-}
\newcommand{\bp}{b^+}
\newcommand{\bm}{b^-}
\newcommand{\si}{\sigma}
\newcommand{\siIj}{\sigma_{I_j}}

\newcommand{\lap}{\la^+}
\newcommand{\lam}{{\la^-}}
\newcommand{\lapm}{{\la^\pm}}
\newcommand{\lapj}{\la^+(j)}
\newcommand{\lamj}{\la^-(j)}
\newcommand{\lapjM}{\la_M^+(j)}
\newcommand{\lamjM}{\la_M^-(j)}
\newcommand{\lapmj}{\la^\pm(j)}
\newcommand{\tlap}{\tilde{\la}^+}
\newcommand{\tlam}{{\tilde{\la}^-}}
\newcommand{\tlapm}{{\tilde{\la}^\pm}}
\newcommand{\tlapj}{\tilde{\la}^+(j)}
\newcommand{\tlamj}{\tilde{\la}^-(j)}
\newcommand{\tlapmj}{\tilde{\la}^\pm(j)}

\newcommand{\SC}{\operatorname{SC}}
\newcommand{\MP}{\operatorname{MP}}
\newcommand{\KMK}{\operatorname{KMK}}

\newcommand{\gap}{\ga^+}
\newcommand{\gam}{\ga^-}
\newcommand{\bgap}{\bar{\ga}^+}
\newcommand{\bgam}{\bar{\ga}^-}
\newcommand{\gapm}{\ga^\pm}
\newcommand{\gapj}{\ga^+(j)}
\newcommand{\gamj}{\ga^-(j)}
\newcommand{\bgapj}{\bar{\ga}^+(j)}
\newcommand{\bgamj}{\bar{\ga}^-(j)}
\newcommand{\gapmj}{\ga^\pm(j)}
\newcommand{\tgap}{\tilde{\ga}^+}
\newcommand{\tgam}{\tilde{\ga}^-}
\newcommand{\tgapm}{\tilde{\ga}^\pm}
\newcommand{\tgapj}{\tilde{\ga}^+(j)}
\newcommand{\tgamj}{\tilde{\ga}^-(j)}
\newcommand{\tgapmj}{\tilde{\ga}^\pm(j)}

\newcommand{\Ep}{E^+}
\newcommand{\Em}{E^-}

\begin{abstract}
In the theory of orthogonal polynomials, sum rules are remarkable 
relationships between a functional defined on a subset of all probability measures 
 involving the reverse Kullback-Leibler divergence with respect to a particular distribution and 
recursion coefficients related to the orthogonal polynomial construction. 
Killip and Simon (\cite{Killip2}) have given a revival interest to this subject by showing a  quite surprising sum rule for measures dominating the semicircular distribution on $[-2,2]$. 
This sum rule includes a contribution of the atomic part of the measure away from $[-2,2]$.
In this paper, we recover this sum rule by using 
probabilistic tools on random matrices. Furthermore,  we obtain new (up to our knowledge) {\it magic} sum rules for the reverse Kullback-Leibler divergence with respect to the Marchenko-Pastur or  Kesten-McKay distributions.  As in the semicircular case, these formulas 
 include a contribution of the atomic part appearing away from the support of the reference measure. 
\smallskip   

\end{abstract}
\maketitle

\section{Introduction}
\label{sectIN}
\subsection{Szeg\H{o}-Verblunsky theorem and sum rules}
\label{susze}
A very  famous result in the theory of orthogonal polynomial on the unit circle (OPUC) is the Szeg\H{o}-Verblunsky theorem (see \cite{Simon-newbook} Theorem 1.8.6 p. 29). It concerns a 
deep relationship between the coefficients involved in the construction of the orthogonal polynomial sequence of a measure supported by the unit circle and its  logarithmic entropy.  More precisely, the inductive relation  between two successive monic orthogonal polynomials $\phi_{n+1}$ and $\phi_{n}$ ($\deg\phi_n=n$, $n\geq 0$) associated with a probability measure $\mu$ on the unit circle $\T$ supported by at least $n+1$ points 
involves a complex number $\alpha_n$
and may be written as
\begin{equation}
\label{recpolycirc}
\phi_{n+1}(z)=z\phi_{n}(z)-\overline{\alpha}_n\phi_{n}^*(z)\mbox{ where } \phi_{n}^*(z):=z^n\overline{\phi_n(1/\bar{z})}.
\end{equation}
The complex number $\alpha_n=-\overline{\phi_{n+1}(0)}$ is the so-called Verblunsky coefficient. In other contexts, it is also called Schur, Levinson, Szeg\H{o} coefficient or even canonical moment (\cite{DeSt97}).

The Szeg\H{o}-Verblunsky  theorem is the identity
\begin{equation}
\label{segverth}
\frac{1}{2\pi}\int_{0}^{2\pi}\log g_{\mu}(\theta)d\theta = \sum_{n\geq 0}\log(1-|\alpha_n|^2)\,,
\end{equation}
where the Lebesgue decomposition 
of $\mu$ with respect to  the uniform measure $d\theta/2\pi$ on $\T$ is
$$d\mu(\theta)= g_{\mu}(\theta)\frac{d\theta}{2\pi}+d\mu_s(\theta)\,,$$
and where both sides of (\ref{segverth}) are simultaneously finite or infinite. 
An exhaustive overview and the genesis tale of this 
crucial theorem may be found in the very nice book of Simon (\cite{Simon-newbook}). The identity (\ref{segverth}) is one of the most representative example of a sum rule (or trace formula): it connects the coefficients of an operator (\cite{Killip3})
to its spectral data. There are various analytical methods of proof (see Chapter 1 in \cite{Simon-newbook}) and a probabilistic one (see section 5.2 of \cite{gamboacanonical}).


In the theory of orthogonal polynomials on the real line (OPRL), given a probability measure $\mu$ with an infinite support, a.k.a. nontrivial case  (resp. with a finite support consisting of $n$ points, a.k.a. trivial case), the orthonormal polynomials  (with positive leading coefficients) obtained by applying the orthonormalizing Gram-Schmidt procedure
to the sequence $1, x, x^2, \dots$ obey the recursion relation
\begin{align} \label{polrecursion}
xp_k(x) = a_{k+1} p_{k+1}(x) + b_{k+1} p_k (x) + a_{k} p_{k-1}(x)
\end{align}
for $ k \geq 0$ (resp. for $0 \leq k \leq n-1$) where the Jacobi parameters satisfy $b_k \in \mathbb R, a_k > 0$. Notice that here the orthogonal polynomials are not monic but normalized in $L^2(\mu)$.


The sum rule analogous to (\ref{segverth}) in the OPRL case is given by the Killip-Simon theorem (\cite{Killip2}).  It relates the sum of functions of $a_k$ and $b_k$ to a spectral expression involving $\mu$. Like the Szeg\H{o}-Verblunsky formula, the spectral side of the sum rule equation measures in some sense the deviation from a reference measure, the semi-circle law
\begin{equation}
\label{SC0}
\operatorname{SC}(dx) = \frac{1}{2\pi}\sqrt{4-x^2}\!\ \mathbbm{1}_{[-2, 2]}(x)\!\ dx
\end{equation} 
and gives on the ``sum-side'' the corresponding contribution by the sequence of recursion coefficients. We restate the sum rule of \cite{Killip2} in Section \ref{susectHER} in full detail. Again, an exhaustive discussion and history of this sum rule can be found in  Section 1.10 of the book of \cite{Simon-newbook}. The deep analytical proof is in Chapter 3 of the book.

In both models, the Kullback-Leibler divergence or relative entropy between two probability measures  
$\mu$ and $\nu$ plays a major role. When the probability space is $\mathbb R$ endowed with its Borel  $\sigma$-field 
it  is defined by
\begin{equation}
\label{KL}
{\mathcal K}(\mu\ |\ \nu)= \begin{cases}  \ \displaystyle\int_{\mathbb R}\log\frac{d\mu}{d\nu}\!\ d\mu\;\;& \mbox{if}\ \mu\ \hbox{is absolutely continuous with respect to}\ \nu\\
   \  \infty  &  \mbox{otherwise.}
\end{cases}
\end{equation}
Usually, $\nu$ is the reference measure. Here  the spectral side  will involve the reversed Kullback-Leibler divergence, where $\mu$ is the reference measure and $\nu$ is the argument. 

Later, \cite{Nazarov} obtained a more general sum rule, when the reference measure is $A(x) \operatorname{SC}(dx)$ with $A$ a nonnegative polynomial (see also \cite{KuKu} for other generalizations). We will discuss this point in Section \ref{conj}.

\subsection{Our main results with hands: outline of the paper}
\label{suhand}

The contribution of this paper 
is twofold. On the one hand, we show two new sum rules. One for measures on the positive half line and one for measures restricted to a compact interval. In each case, the reference measure is different and the sum involves
a function of specific coefficients related to the sequences $(a_n)_n$ and $(b_n)_n$. On the other hand, we also show new large deviation theorems for spectral measures of random operators. In fact, this probabilistic result yields the new sum rules as a direct consequence and also allows for an alternative probabilistic proof of the Killip-Simon sum rule. Notice that as pointed out by Simon
in \cite{gem}:
{\it "The gems of spectral theory
are ones
that set up one-one correspondences between classes of measures and
coefficients with some properties."} In Section \ref{sgemgem}, we will discuss the underlying gems deduced from the new sum rules.  

Large deviations for these random spectral measures arising from the classical ensembles of random matrix theory have been considered before in \cite{FGAR}. Therein, the main tool was the study of large deviation properties of the recursion coefficients. This method yields as rate function precisely the {\it sum side} of the new sum rules. Furthermore, the rate function in our new result 
is the {\it spectral side}. Since the rate function in large deviations is unique, both sides must be equal. This is, in a nutshell, our 
 proof for sum rules. Our method of proof also 
stress why both sides of the sum rule equations gives a measure for the 
{\it divergence} to some reference measure as they are large deviation rate functions.

To build a comfortable common ground for both mathematical analyst and probabilist reader, we will recap in the further course of this section some useful facts on spectral measures, on their tridiagonal representations and on their randomization. 
We will also recall the definition of large deviations. In Section \ref{sectSU}, we restate the sum rule showed in \cite{Killip2} and give our new sum rules. For the convenience of the reader, we formulate the sum rules without mentioning the underlying randomization. In Section \ref{sectLAR} we give the main large deviation result and we explain why the sum rules are a consequence of this theorem. We also give a conjecture for a more general sum rule going away from the frame of classical ensembles. The proof of the main large deviation theorem can be found in Section \ref{sectproof}. Finally, some technical details are referred to the Appendix.

Let us notice that three extensions of the above method are quite natural and will appear in further work. 
\begin{enumerate}
\item A matricial version of the Killip-Simon sum rule
is the due to Damanik, Killip and Simon (see \cite{damal} or Theorem 4.6.3 in \cite{Simon-newbook}). We will extend the results of the present paper to block Laguerre and block Jacobi  random matrices. As a matter of fact, we will lean on the large deviation results  proved in \cite{JMVA}.
\item In the unit circle case, there is a natural model having a limit measure supported by a proper arc of $\mathbb T$. In this frame, the random Verblunsky coefficients have a nice independence structure (see \cite{BNR}). This allows to extend the sum rules developed here. 
\item All along this paper, we consider measures with essential support  consisting in a single interval (so-called one-cut assumption). We will later consider equilibrium measures supported by several intervals. For this task, the probabilistic tools may be found in \cite{borot-multi} and the analytic ones are in Chapter 9 of \cite{Simon-newbook}.
\end{enumerate}
\subsection{OPRL and tridiagonal matrices}


If $H$ is a self-adjoint bounded operator on a Hilbert space $\mathcal H$ and $e$ is a cyclic vector, the spectral measure of the  pair $(H,e)$ is the unique probability measure $\mu$ on $\mathbb R$ such that
\[\langle e, H^k e\rangle = \int_\mathbb R x^k d\mu(x) \ \ (k \geq 1).\]
Actually, $\mu$ is a unitary invariant for $(H,e)$. Another invariant is the tridiagonal reduction whose coefficients will play the role of the earlier-mentioned Verblunsky coefficients for unitary operators. If dim $\mathcal H =n$ and $e$ is cyclic for $H$,  let $\lambda_1, \ldots, \lambda_n$ be the (real) eigenvalues of $H$ and let $\psi_1, \ldots, \psi_n$ be a system of orthonormal eigenvectors. The spectral measure of the pair $(H,e)$ is then
\begin{align}\label{spectralmeasure}
\mun =  \sum_{k=1}^n \w_k\delta_{\lambda_k}\,,
\end{align}
with $\w_k= |\langle \psi_k, e\rangle|^2$. This measure is a weighted version of the empirical eigenvalue distribution
\begin{align}\label{empiricallaw}
\muun = \frac{1}{n} \sum_{k=1}^n \delta_{\lambda_k} \,.
\end{align}
$\mun$ is called {\it eigenvector empirical distribution function} in a recent paper of  \cite{XQB}. 
Let us now describe shortly the Jacobi mapping between tridiagonal matrices and spectral measures.

 We consider $n \times n$ 
 matrices  corresponding to  measures supported by 
$n$ points (trivial case) and semi-infinite matrices corresponding to measures with bounded infinite support (non-trivial case). 
 
 In the basis $\{p_0, p_1, \dots, p_{n-1}\}$, the 
linear transform $f(x) \rightarrow xf(x)$ (multiplication by $x$) 
 in $L^2(d\mu)$ is represented by the matrix
\ben
\label{favardfini}
J_\mu = \begin{pmatrix} b_1&a_1 &0&\dots&0\\
a_1&b_2  &a_2&\ddots&\vdots\\
0&\ddots &\ddots&\ddots&0\\
\vdots&\ddots&a_{n-2}&b_{n-1}&a_{n-1}\\
0&\dots&0&a_{n-1}&b_{n}
\end{pmatrix}
\een
So, measures supported by $n$ points lead to Jacobi matrices, i.e. $n\times n$ symmetric tridiagonal matrices with subdiagonal positive terms. In fact, there
is a one-to-one correspondence between such a matrix and such a measure. 
If $J$ is  a Jacobi matrix, 
we can take the first vector of the canonical basis as the cyclic vector $e$. 
 Let  $\mu$  be the spectral measure associated to the pair $(J, e_1)$, then $J$ represents the multiplication by $x$ in the basis of orthonormal polynomials associated to $\mu$ and $J=J_\mu$.

More generally, if $\mu$ is a probability measure on $\R$ with bounded infinite support, we may apply the same Gram-Schmidt process and consider the associated semi-infinite Jacobi matrix: 
\ben
\label{favardinfini}
J_\mu = \begin{pmatrix} b_1&a_1 &0&0&\cdots\\
a_1&b_2  &a_2&0&\cdots\\
0&a_2 &b_3&a_3& \\
\vdots& &\ddots&\ddots&\ddots
\end{pmatrix}
\een
Notice that  again we have $a_k > 0$ for every $k$. The mapping $\mu \mapsto J_\mu$ (called here the Jacobi mapping) 
 is a one to one correspondence between probability measures on $\mathbb R$ having compact infinite support and this kind of tridiagonal matrices 
with  $\sup_n(|a_n| +
|b_n|) < \infty$. This result is 
sometimes called Favard's theorem. 

\subsection{Randomization: gas distribution and random matrices}   
\label{suserand}

In this paper we consider distributions of log-gases and random matrices. In the sequel, $n$ is the number of particles (or eigenvalues), denoted by $\lambda_1,\dots ,\lambda_n$, with the joint distributions $\Pnv$ on $\R^n$ having  the 
 density
\begin{align}\label{generaldensity}
\frac{d\Pnv(\lambda)}{d\lambda} = \frac{1}{Z_V^n}
e^{- n\beta' \sum_{k=1}^nV(\lambda_k)}\prod_{1\leq  i < j\leq n} |\lambda_i - \lambda_j|^\beta.
\end{align}
with respect to the Lebesgue measure $d\lambda = d\lambda_1 \cdots d\lambda_n$. 
The potential $V:\mathbb{R} \to (-\infty,+\infty]$ is supposed to be continuous real valued on the interval ${(b^-, b^+)}$ ($-\infty\leq b^-<b^+\leq+\infty$), infinite outside of $[b^-,b^+]$ and $\lim_{x\to b^\pm} V(x) = V(b^\pm)$ with possible limit $V(b^\pm)=+\infty$. Let $\beta=2\beta'>0$ be the inverse temperature. Under the assumption 
\begin{itemize}
\item[(A1)] Confinement: 
$\displaystyle \qquad
\liminf_{x \rightarrow b^\pm} \frac{V(x)}{2 \log |x|} > \max(1, \beta^{-1}) \, , $
\end{itemize}
the empirical distribution $\muun$ of eigenvalues $\lambda_1,\dots ,\lambda_n$ has a 
 limit $\mu_V$ (in probability)\footnote{Various authors used to say  {\it almost surely}, but since the probability spaces are not embedded, it seems more convenient to keep {\it in probability}.}, which is the unique minimizer of 
\begin{align}
\label{ratemuu}
\mu  \mapsto \mathcal E (\mu) := \int V(x) d\mu(x) - \iint \log |x-y| d\mu(x)d\mu(y).
\end{align}
$\mu_V$ has compact support (see \cite{johansson1998fluctuations} or \cite{agz}). 
Indeed, this is a consequence of the large deviations of the empirical spectral measure (see Theorem \ref{LDPmuu}). 
We will make the following assumptions on $\mu_V$:
\begin{itemize}
\item[(A2)] One-cut regime: the support of $\mu_V$ is a single interval $[\alpha^-, \alpha^+]\subset [b^-, b^+]$ ( $\alpha^-< \alpha^+$).
\item[(A3)] Control (of large deviations): the effective potential
\begin{align}
\label{poteff}
\Ir_V (x) := V(x) -2\int \log |x-\xi|\!\ d\mu_V(\xi)
\end{align}
achieves its global minimum value on $(b^-, b^+) \setminus (\alpha^-, \alpha^+)$  only on the boundary of this set.
\end{itemize}
Furthermore, to obtain a non-variational expression for the rate we need the following conditions:
\begin{itemize}
\item[(A4)] Offcriticality: We have
\begin{align*}
d\mu_V(x) = \frac{1}{2\pi}S(x) \sqrt{ 
\frac{\Pi_{\tau \in \operatorname{Soft}} |x-\alpha^\tau | }{\Pi_{\tau' \in \operatorname{Hard}} |x-\alpha^{\tau'} |}
 }\!\ dx
\end{align*}
where $S> 0$ on $[\alpha^-, \alpha^+]$ and $\tau \in \operatorname{Hard}$ iff $b^\tau=\alpha^\tau$, otherwise $\tau \in \operatorname{Soft}$ ($\operatorname{Hard}\cap \operatorname{Soft}=\emptyset$ and $\operatorname{Hard}\cup \operatorname{Soft}=\{-,+\}$).
\item[(A5)] Analyticity: $V$ can be extended as an holomorphic function is some open neighborhood of $[\alpha^-, \alpha^+]$.
\end{itemize}
We remark that for $V$ strictly convex, the assumptions (A2), (A3) and (A4) are fulfilled (see \cite{Borot} and \cite{johansson1998fluctuations}).

Hereafter, we discuss the classical models with their potentials  and their domains, and 
 their equilibrium measure.
\begin{enumerate}
\item Hermite ensemble:
$$V(x) = \frac{x^2}{2} \ \ , \ \  (b^-, b^+)= \R \mbox{ and the equilibrium measure is } \hbox{SC}(dx)  , \ \alpha^\pm = \pm 2\,.$$
\item  Laguerre ensemble of parameter $\tau \in (0,1]$:
\[V(x) = \tau^{-1} x - (\tau^{-1} -1) \log x  \ \ , \ \  [b^-, b^+) = [0, \infty)\]
with equilibrium measure the Marchenko-Pastur law with parameter $\tau$,
\[\hbox{MP}_\tau(dx) = \frac{\sqrt{(\tau^+ -x)(x-\tau^-)}}
{2\pi \tau x} \ \mathbbm{1}_{(\tau^-, \tau^+)} (x) dx\ \ , \ \ \alpha^\pm=\tau^\pm = (1\pm \sqrt\tau)^2\]
\item Jacobi ensemble of parameters $\kappa_1, \kappa_2\geq 0$:
\[V(x) = -\kappa_1 \log x - \kappa_2 \log (1-x)  \ \ , \ \  [b^-, b^+] = [0,1]\,.\]
The equilibrium measure is the Kesten-McKay distribution 
\[\hbox{KMK}_{\kappa_1, \kappa_2}(dx) = \frac{(2+\kappa_1+\kappa_2)}{2\pi}\frac{\sqrt{(u^+ -x)(x- u^-)}}{x(1-x)} \ \mathbbm{1}_{(u^-, u^+)}(x) dx\, ,\]
where
\begin{align}
\label{upm}
\alpha^\pm= u^\pm := \frac{1}{2} + \frac{\kappa_1^2 - \kappa_2^2 \pm 4 \sqrt{(1+\kappa_1)(1+\kappa_2)(1+\kappa_1+\kappa_2)}}{2(2+\kappa_1+\kappa_2)^2} .
\end{align}
\end{enumerate}

Let us start with a crash recall in random matrix theory. The GOE  of order $n$ is a probability distribution $P_n^{(1)}$ on the set of all symmetric real $n \times n$ matrices, obtained by assuming that the diagonal entries are distributed as $\mathcal N (0,2)$ and the non-diagonal ones as $\mathcal N (0,1)$ and that entries are independent up to symmetry. Taking on-or-above-diagonal entries as coordinates of the random matrix $H$, this gives a density with respect to the Lebesgue measure proportional to
 $\exp - \tr H^2/4$.  
The distribution of eigenvalues of $\frac{1}{\sqrt{n}} H$ is given by (\ref{generaldensity}) with $\beta=1$ and $V(x) = x^2/2$. Besides, 
\cite{Trotter} proved that the  coefficients $b_k$ are Gaussian and the coefficients $a_k^2$ are gamma distributed, with convenient parameter. Furthermore, by the invariance of the Gaussian distribution under rotation, the  first row of the eigenvector matrix is independent of the eigenvalues and uniformly distributed on the sphere. \cite{Dumitriu1} also proved that conversely, if we take an array $(\lambda_1, \dots, \lambda_n)$ distributed as in (\ref{generaldensity}), with general $\beta > 0$ and an independent array of weights $(\w_1, \dots, \w_n)$ sampled with the Dirichlet  distribution Dir$_n(\beta')$ of order $n$ and  parameter $\beta'$ 
 on the simplex $\sum_i \w_i = 1$, i.e. with density proportional to 
\[\left( \w_1 \cdots \w_n \right)^{\beta' -1},\]
then the coefficients of the tridiagonal matrix are independent  Gaussian and gamma variables, respectively. In the case of Laguerre and Jacobi ensembles, other systems of auxiliary variables with nice structure of independence were introduced in  \cite{Dumitriu1} and in \cite{Killip1}. We will use these parametrizations in the next sections.

The case of a general potential $V$ is not so easy. Nevertheless, the correspondence is ruled by the following result. Let $A\sn = (a_1, \dots, a_{n-1})$ and $B\sn = (b_1, \dots, b_n)$ and let
 $T\sn$ 
 be the symmetric tridiagonal matrix with $T^{(n)}_{k,k} =b_k$ for $k \leq n$ and $T^{(n)}_{k, k+1} =T\sn_{k+1,k} = a_{k}$ for $k\leq n-1$. Then we may state the following theorem.

\begin{thm}[\cite{krishna13} Prop.2]
\label{KRV}
Let $(A\sn,B\sn)$ sampled from the density proportional to
\[\exp -n \beta' \left[\tr\  V(T\sn) - 2\sum_{k=1}^{n-1} \left( 1 - \frac{k}{n} - \frac{1}{n\beta}\right) \log a_k \right]\,.\] 
Then the eigenvalues $\lambda_1, \dots, \lambda_n$ have joint density proportional to
\[e^{- n \beta' 
\sum_1^n V(\lambda_k)}\ \prod_{1\leq j <k\leq n}|\lambda_j - \lambda_k|^\beta \]
and the weights $\w_1, \dots, \w_n$ are independent with distribution $\Dir_n(\beta')$.
\end{thm}

\subsection{Large deviations}

\label{sular}

In order to be self-contained, let us recall the definition of a large deviation principle. For a general reference of large deviation statements we refer to the book of \cite{demboz98} or to the Appendix D of \cite{agz}.

 Let $U$ be a topological Hausdorff space with Borel $\sigma$-algebra $\mathcal{B}(U)$. We say that a sequence $(P_{n})_n$ of probability measures on $(U,\mathcal{B}(U))$ satisfies a large deviation principle (LDP) with speed $a_n$ and
rate function $\mathcal{I} : U \rightarrow [0, \infty]$  if:
\begin{itemize}
\item [(i)] $\mathcal I$ is lower semicontinuous.
\item[(ii)] For all closed sets $F \subset U$:
\begin{align*}
\limsup_{n\rightarrow\infty} \frac{1}{a_n} \log P_{n}(F)\leq -\inf_{x\in F}\mathcal{I}(x)
\end{align*}
\item[(iii)] For all open sets $O \subset U$:
\begin{align*}
\liminf_{n\rightarrow\infty} \frac{1}{a_n} \log P_{n}(O)\geq -\inf_{x\in O}\mathcal{I}(x)
\end{align*}
\end{itemize}
The rate function $\mathcal{I}$ is good if its level sets
$\{x\in U |\ \mathcal{I}(x)\leq a\}$ are compact for all $a\geq 0$. 
If in the conditions above, we replace {\it closed sets} by {\it compact sets}, we say that $(P_{n})_n$ satisfies a weak LDP. In this case, we can recover a LDP if the additional condition of exponential tighness is fulfilled:

For every $M > 0$ there exists a compact set $K_M \subset U$ such that
\[\limsup_{n\rightarrow\infty}  \frac{1}{a_n} \log P_{n}(U \setminus K_M) \leq -M\,.\]
In our case, the measures $P_n$ will be the distributions of the random spectral measures $\mu_n$ and we will say that the sequence of measures $\mu_n$ satisfies a LDP. 

The most famous LDP in random matrix theory  is for the sequence of empirical spectral measures. We let $\mathcal{P}_1$ denote the set of all probability measures on $\mathbb{R}$.
\begin{thm}
\label{LDPmuu}
If the potential $V$ satisfies assumption (A1), and if $(\lambda_1, \dots, \lambda_n)$ is distributed according to $\Pnv$ (see (\ref{generaldensity})),
then the sequence of random probability measures $(\mu_\u\sn)$
satisfies in $\mathcal{P}_1$ equipped with the weak topology, a LDP with speed $\beta' n^2$ and good rate function
\begin{align*}\mu \mapsto \mathcal E (\mu) - \mathcal E (\mu_V)\end{align*} where $\mathcal E$  is defined  in (\ref{ratemuu})\,. \end{thm}
Let us recall the definition of convex duality, used several times in the proofs. If $R$ is a function defined on a topological vector space $\mathcal H$ and valued in $(-\infty, \infty]$, then its convex dual $R^*$  is a function defined on the topological dual space $\mathcal H^*$ by
\[R^* (x) = \sup_{\theta \in \mathcal H} \left[\langle \theta , x\rangle - R(\theta)\right]\quad (x\in\mathcal H^*).\]
Here, $\langle\ , \ \rangle$ is the duality bracket. For $\mathcal H=\R$,  two examples are meaningful in our context
\begin{enumerate}
\item Gaussian case
\begin{equation}
\label{Laplace1}
L_0(\theta) = \frac{\theta^2}{2} \Rightarrow L_0^*(x) = \frac{x^2}{2} \end{equation}
\item Exponential case
\begin{equation}
\label{LaplaceL}
 L(\theta) = \begin{cases} \  -\log(1-\theta)& \mbox{ if}\ \theta < 1\,, \\
\ \infty & \mbox{ otherwise,}
\end{cases}
\end{equation}
then
\begin{equation}
\label{Laplace2}
L^*(x) = G(x) := \begin{cases} \ x- 1 - \log x& \mbox{ if}\  x >0\,, \\
\ \infty & \mbox{ otherwise.}
\end{cases}
\end{equation}
\end{enumerate}
Let us denote respectively by $\mathcal N(0, \sigma^2)$ the centered Gaussian distribution and Gamma$(a,b)$ the gamma distribution of order $a> 0$ and scale factor $b > 0$  with density
\[x \mapsto \frac{x^{a-1}}{b^a \Gamma(a)} e^{-\frac{x}{b}} \ \  x > 0\, .\] 
$L_0^*$ 
 is the rate function of the LDP satisfied by  $(\mathcal N(0, n^{-1}))_n$ at speed $n$
 and $L^*$  
 is the rate function of the LDP satisfied by  (Gamma$(n, n^{-1}))_n$, also at speed $n$.
Besides, (Gamma$(a, n^{-1}))_n$ satisfies a LDP at speed $n$ with rate function 
\begin{equation}
\label{LDPdeg}
x \mapsto  \begin{cases}\  x& \mbox{ if}\  x \geq 0 , \\
\ \infty & \mbox{ otherwise.}
\end{cases}
\end{equation}

\section{Sum rules from large deviations}
\label{sectSU}

Let $\Sr = \Sr(\am,\ap)$ be the set of all bounded positive measures $\mu$ on $\R$ with 
\begin{itemize}
\item[(i)] $\operatorname{supp}(\mu) = J \cup \{\lambda_i^-\}_{i=1}^{N^-} \cup \{\lambda_i^+\}_{i=1}^{N^+}$, where $J\subset I= [\am,\ap]$, $N^-,N^+\in\N\cup\{\infty\}$ and 
\begin{align*}
\lambda_1^-<\lambda_2^-<\dots <\am \quad \text{and} \quad \lambda_1^+>\lambda_2^+>\dots >\ap .
\end{align*}
\item[(ii)] If $N^-$ (resp. $N^+$) is infinite, then $\lambda_j^-$ converges towards $\am$ (resp. $\lambda_j^+$ converges to $\ap$).
\end{itemize}
Such a measure $\mu$ will be written as
\begin{align}\label{muinS0}
\mu = \mu_{|I} +  \sum_{i=1}^{N^+} \gamma_i^+ \delta_{\lambda_i^+} + \sum_{i=1}^{N^-} \gamma_i^- \delta_{\lambda_i^-}
\end{align}
Further, we define $\Sr_1=\Sr_1(\am,\ap):=\{\mu \in \Sr |\, \mu(\R)=1\}$. We endow $\Sr_1$ with the weak topology and the corresponding Borel $\sigma$-algebra.

\subsection{Hermite case revisited}
\label{susectHER}

\medskip

We start by stating the classical sum rule (due to  \cite{Killip2} and explained in \cite{Simon-newbook} p.37), the new probabilistic proof using large deviations is tackled  in Section \ref{sectLAR2}. The sum rule gives two different expressions for the {\it distance} to the semicircle law $\operatorname{SC}$. Its Jacobi coefficients are
\begin{equation}
\label{JcH}a^{\SC}_k = 1,\  b^{\SC}_k =0 \ \ \hbox{for all} \ k\geq 1\,.\end{equation}
For a probability measure $\mu$ on $\R$ with recursion coefficients $(a_k)_k,\, (b_k)_k$ as in \eqref{polrecursion}, define the sum
\begin{align}\label{rateG1}
\mathcal{I}_H(\mu) = \sum_{k\geq 1}\big( \frac{1}{2} b_k^2 + G(a_k^2)\big)= \sum_{k\geq 1} \big(L_0^*(b_k) + G(a_k^2)\big),
\end{align}
where $G$ and $L_0^*$ have been defined in the previous section. 
Further,  let
\begin{align*}
\mathcal{F}_H^+(x) :=  \begin{cases} &  \displaystyle \int_2^x \sqrt{t^2-4}\!\ dt = \tfrac{x}{2} \sqrt{x^2-4} - 2 \log \left( \tfrac{x+\sqrt{x^2-4}}{2}\right)\;\;\;\;\mbox{if} \ x \geq 2\\
    &  \infty\;\;\mbox{ otherwise,}
\end{cases}
\end{align*}
and, for $x\in\R$,  set $\mathcal{F}_H^-(x):=\mathcal{F}_H^+(-x)$. Then we have the following theorem.

\medskip

\begin{thm}[\cite{Killip2}]
\label{sumruleg}
Let $J$ be a Jacobi matrix with diagonal entries $b_1,b_2,\ldots \in \R$ and subdiagonal entries $a_1,a_2,\ldots >0$ satisfying  $\sup_k a_k + \sup_k |b_k| < \infty$ and let $\mu$ be the associated spectral measure. Then $\mathcal{I}_H(\mu)$ is infinite if $\mu \notin \Sr_1(-2,2)$ and for $\mu \in \Sr_1(-2,2)$, 
\begin{align*} 
\mathcal{I}_H(\mu) = {\mathcal K}(\operatorname{SC}\! |\!\ \mu) +  \sum_{n=1}^{N^+} {\mathcal F}^+_H(\lambda_n^+)  +  \sum_{n=1}^{N^-} {\mathcal F}^-_H(\lambda_n^-)
\end{align*}
where both sides may be infinite simultaneously. 
\end{thm}

\medskip

\subsection{New magic sum rule: the Laguerre case}
\label{susectLAG}

For our first new sum rule, we consider probability measures $\mu$ supported on $[0,\infty)$. In this case, the recursion coefficients can be decomposed as
\begin{align} \label{zerl1}
\begin{split}
b_k =& z_{2k-2} + z_{2k-1}, \\
a_k^2 =&  z_{2k-1}z_{2k},
\end{split}
\end{align}
for $k \geq 1$, where $z_k\geq 0$ and $z_{0}=0$. In fact, by Favard's Theorem a measure $\mu$ is supported on $[0,\infty)$ if and only if the decomposition as in \eqref{zerl1} holds. The central probability measure is the Marchenko-Pastur law $\operatorname{MP}_\tau$ defined in Section \ref{suserand}, whose Jacobi coefficients  are
\begin{equation}
\label{JcL}a^{\MP}_k = \sqrt \tau \ (k\geq 1) \ \ , \ \ b^{\MP}_1 = 1 \ , \ b^{\MP}_k = 1+\tau \ \ (k \geq 2)\end{equation}
and correspond to $z^{\MP}_{2n-1}=1$ and 
$z^{\MP}_{2n} = \tau$ for all $n \geq 1$. 
For a measure $\mu$ supported on $[0,\infty)$, let
\begin{align} \label{rateL1}
\mathcal{I}_L(\mu) :=   \sum_{k=1}^\infty \tau^{-1}G(z_{2k-1}) +  G(\tau^{-1}z_{2k}) . 
\end{align}
For the new sum rule, we have to replace $\mathcal{F}_H^\pm$ by
\begin{align*}
\mathcal{F}_L^+(x) =  \begin{cases} &  \displaystyle \int_{\tau^+}^x \frac{\sqrt{(t -  \tau^-)(t - \tau^+)}}{t\tau}\!\ dt \;\;\;\;\mbox{if} \ x \geq \tau^+,\\
    &  \infty\;\;\mbox{ otherwise,}
\end{cases}
\end{align*}
\begin{align*} 
{\mathcal F}_L^-(x) =  \begin{cases} & \displaystyle\int_x^{\tau^-} \frac{\sqrt{(\tau^- -t)(\tau^+ -t)}}{t\tau}\!\ dt  \;\;\;\;\mbox{if} \ x \leq \tau^-,\\
    &  \infty\;\;\mbox{ otherwise.}
\end{cases}  
\end{align*}
Then we have the following magic sum rule for probability measures on $[0,\infty)$. The 
 probabilistic proof can be found in Section \ref{sectLAR2}.

\medskip

\begin{thm}\label{sumrulel}
Assume the entries of the Jacobi matrix $J$ can be decomposed as in \eqref{zerl1} with $\sup_k z_k <\infty$ and let $\mu$ be the spectral measure of $J$. Then for all $\tau \in (0,1]$, $\mathcal{I}_L(\mu)=\infty$ if $\mu \notin \Sr_1(\tau^-,\tau^+)$. If $\mu \in \Sr_1(\tau^-,\tau^+)$, we have 
\begin{align*} 
\mathcal{I}_L(\mu) = {\mathcal K}(\operatorname{MP}
_{\tau}
\! |\!\ \mu) +  \sum_{n=1}^{N^+} {\mathcal F}_L^+(\lambda_n^+)  +  \sum_{n=1}^{N^-} {\mathcal F}_L^-(\lambda_n^-)
\end{align*}
where both sides may be infinite simultaneously. 
\end{thm}

\medskip

Note that if $\tau=1$, the support of the limit measure is $[0,4]$, so that we have a hard edge at 0 with $N^-=0$ and no contribution of outliers to the left.

\medskip

\subsection{New magic sum rule: the Jacobi case}
\label{susectJAC}

Our second new sum rule is a generalization of the Szeg\H{o} theorem for probability measures on the unit circle. 
The classical Szeg\H{o} mapping is a correspondence between a probability measure $\nu$ on $\mathbb T$ invariant by $\theta \mapsto 2\pi -\theta$ and a probability measure $\mu$ on $[-2, 2]$ obtained by pushing forward $\nu$ by 
 the mapping $\theta \mapsto 2 \cos \theta$. In this case the Verblunsky coefficients $(\alpha_k)_{k\geq 0}$ of $\nu$ (they all belong to $[-1, 1]$ by symmetry) are by extension called the Verblunsky coefficients of $\mu$. For $k\geq 1$, the recursion coefficients associated with $\mu$ are connected with the Verblunsky coefficients  by the Geronimus relations: 
%
\begin{align}\label{zerl2}
\begin{split}
b_{k+1} &= (1-\alpha_{2k-1})\alpha_{2k}-(1+\alpha_{2k-1})\alpha_{2k-2}\\
a_{k+1} &= \sqrt{(1-\alpha_{2k-1})(1-\alpha_{2k}^2)(1+\alpha_{2k+1})}
\end{split}
\end{align}
where $\alpha_k \in [-1,1]$ and $\alpha_{-1}=-1$. While these recursion coefficients give a measure $\mu$ on $[-2,2]$, it is more convenient for our approach to consider the measure $\tilde\mu$ on $[0,1]$ obtained by  pushing forward  $\mu$ 
 by the affine mapping $x\mapsto \tfrac{1}{2}-\tfrac{1}{4}x$. 
We keep calling  $(\alpha_k)_k$ the Verblunsky coeffcients of $\tilde \mu$. The Jacobi coefficients of $\widetilde \mu$ are 
\[\widetilde b_k  = \frac{2-b_k}{4} \ , \ \widetilde a_k = \frac{a_k}{4} \ \ (k \geq 1)\,.\]

Here, the important probability measure is the Kesten-McKay distribution $\operatorname{KMK}_{\kappa_1,\kappa_2}$ on $[0,1]$ with parameters $\kappa_1,\kappa_2\geq 0$.
The associated Verblunsky coefficients 
 are, for $k\geq 0$,
\begin{align*}
\alpha_{2k}^{KMK} = \frac{\kappa_1-\kappa_2}{2+\kappa_1+\kappa_2}, \quad \alpha_{2k+1}^{KMK} = -\frac{\kappa_1+\kappa_2}{2+\kappa_1+\kappa_2}.
\end{align*}
and the corresponding Jacobi coefficients are
\begin{equation}
\label{JcJ}\widetilde a^{KMK}_1 =  \frac{\sqrt{(1+\kappa_1)(1+\kappa_2)}}{(2+\kappa_1+\kappa_2)^{3/2}}\ ,\quad \widetilde b^{KMK}_1 =  \frac{1+\kappa_2}{2+\kappa_1+\kappa_2}\,,\end{equation}
and for $k \geq 2$
\[\widetilde a^{KMK}_k = \frac{\sqrt{(1+\kappa_1 + \kappa_2)(1+\kappa_1)(1+\kappa_2)}}{(2+ \kappa_1+ \kappa_2)^2}\ ,\quad    \widetilde b^{KMK}_k=\frac{1}{2}\left[1 -  \frac{\kappa_1^2- \kappa_2^2}{(2+ \kappa_1+ \kappa_2)^2}\right]\,.\]
Set
\begin{align}\label{rateJ1}
\mathcal{I}_J(\tilde\mu) =  \sum_{k=0}^\infty H_1(\alpha_{2k+1}) + H_2(\alpha_{2k}) ,
\end{align}
where for $x\in [-1,1]$
\begin{align*}
H_1(x) = -(1+\kappa_1+\kappa_2) \log \left[\frac{2+\kappa_1+\kappa_2}{2(1+\kappa_1+\kappa_2)} (1-x)\right]
\ -\log \left[\frac{2+\kappa_1+\kappa_2}{2}(1+x)\right]
ù\end{align*}
\begin{align*}
H_2(x) =  - (1+\kappa_1)\log \left[\frac{(2 + \kappa_1 + \kappa_2)}{2(1+\kappa_1)}(1+x)\right]
\  - (1+ \kappa_2)
 \log\left[  \frac{(2 + \kappa_1 + \kappa_2)}{2(1+\kappa_1)}(1-x)\right]\,. 
\end{align*}
Let ${\mathcal F}_J^+$ be defined by
\begin{align*}
{\mathcal F}_J^+(x) = \begin{cases} \ \displaystyle \int_{u^+}^x \frac{\sqrt{(t -  u^+)(t - u^-)}}
{t(1-t)}\!\ dt  & \mbox{ if} \  u^+ \leq  x \leq 1 \\
\ \infty & \mbox{ otherwise.}
\end{cases}
\end{align*}
 Similarly, let
\begin{align*}
{\mathcal F}_J^-(x) = \begin{cases}\ \displaystyle \int_x^{u^-} \frac{\sqrt{(u^--t)(u^+ -t)}}
{t(1-t)}\!\ dt & \mbox{ if} \  0 \leq x \leq u^-\\
\ \infty & \mbox{ otherwise.}
\end{cases}
\end{align*}
Then the following magic sum rule for probability measures on $[0,1]$ holds.

\medskip

\begin{thm}\label{sumrulej}
Let $\tilde\mu$ be a probability measure on $[0,1]$ and let $\alpha_k$ be the Verblunsky coefficients of $\tilde\mu$. Then for any $\kappa_1,\kappa_2\geq 0$, $\mathcal{I}_J(\tilde{\mu})=\infty$ if $\tilde\mu\notin \Sr_1(u^-,u^+)$. If $\tilde\mu\in\Sr_1(u^-,u^+)$, then
\begin{align*} 
\mathcal{I}_J(\tilde\mu) = {\mathcal K}(\operatorname{KMK}_{\kappa_1,\kappa_2}\! |\!\ \tilde\mu) +  \sum_{n=1}^{N^+} {\mathcal F}_J^+(\lambda_n^+)  +  \sum_{n=1}^{N^-} {\mathcal F}_J^-(\lambda_n^-)
\end{align*}
and both sides may be infinite simultaneously. 
\end{thm}

\medskip

Similar to the Laguerre case, if $\kappa_1=0$ or $\kappa_2=0$, then $u^-=0$ or $u^+=1$, respectively, and we have no contribution coming from respective outliers. In particular, if $\kappa_1=\kappa_2=0$, the Kesten-McKay distribution reduces to the arcsine distribution 
\[d\mu_0 (x) = \frac{1}{\pi\sqrt{x(1-x)}} \mathbbm{1}_{(0,1)}(x)\!\ dx\] and then the sum rule reads
\begin{equation}
{\mathcal K}(\mu_0\!\ |\!\ \tilde\mu)  = - \sum_{n=0}^\infty \log (1-\alpha_n^2) ,
\label{sese14}
\end{equation}
which is nothing more than the classical Szeg\H{o} sum rule written for  probability measures pushed forward
by the application $\theta\rightarrow\frac{1}{2}-\frac{\cos\theta}{2}$.   
Notice also that in this frame we may rewrite this sum rule 
in terms of a cousin parametrization. Namely, by the way of the so-called canonical moments defined for a measure supported on $[0,1]$ (see the excellent book of 
 \cite{DeSt97}).   More precisely, let for $k\geq 1$, $p_k$ denote the canonical moment of order $k$ of $\tilde{\mu}$. Recall that $p_k$ may be defined browsing different
paths, the straightest is from the ordinary moments. Indeed, assuming that $\tilde{\mu}$ is not supported by a finite number of points, we have
\begin{align*}
p_1:=&\int_0^1 x\tilde{\mu}(dx), \\
p_{n+1}:=&\frac{\int_0^1 x^{n+1}\tilde{\mu}(dx)-c_{n+1}^-}{c_{n+1}^+-c_{n+1}^-} \quad  \mbox{ for } n\geq 1.
\end{align*}
Here, $c_{n+1}^+$ (resp. $c_{n+1}^-$) is the maximum (resp. minimum) possible value for the $(n+1)$-th moment of a probability measure supported by $[0,1]$ having the same $n$ first moments as $\tilde{\mu}$. With this parametrisation, as $\alpha_n=2p_{n+1}-1$  for all integer $n$ (see \cite{DeSt97} p. 287), the sum rule (\ref{sese14}) becomes

$$ {\mathcal K}(\mu_0\! \ |\!\ \tilde\mu)= 
- \sum_{n=0}^\infty \log (4p_n(1-p_n)).
$$
Notice that the last expression is also the functional obtained in \cite{Gamboa} in the study of large deviations for 
a random Hausdorff moment problem.

\subsection{Semiprecious gems}
\label{sgemgem}
In the introduction of the paper we pointed out that the gem of spectral theory is to set up one-one correspondences between classes of measures and coefficients with some properties. More precisely, 
a gem (see \cite{Simon-newbook} Section 1.4) gives a one-one correspondence between properties on the sequences $(a_n)$ and $(b_n)$ and the associated spectral measure. The gem corresponding to the Hermite case has been proved by Killip and Simon (see \cite{Simon-newbook} Section 1.10). 
We discuss here such correspondences for the two sum rules given in Theorems
\ref{sumrulel} and \ref{sumrulej}. Although we do not succeed to find the holy grail of such a correspondence between classes relying on $(a_n)$ and $(b_n)$,  we set it in terms of the sequences $(z_n)$ (Laguerre case) and $(\alpha_n)$ (Jacobi case).

\begin{cor} \label{semigemL}
Assume the entries of the Jacobi matrix $J$ can be decomposed as in \eqref{zerl1} with $\sup_k z_k <\infty$ and let $\mu$ be the spectral measure of $J$. Then
\begin{align}
\label{zl2}
\sum_{k=1}^\infty[(z_{2k-1}-1)^2 + (z_{2k} - \tau)^2] < \infty
\end{align}
(that is, $\mathcal{I}_L(\mu)< \infty$) if and only if
\begin{enumerate}
\item 
  $\mu \in \mathcal S_1 (\tau^-, \tau^+)$
\item $\sum_{i=1}^{N^+} (\lambda_i^+ - \tau^+)^{3/2} + \sum_{i=1}^{N^-} (\tau^- - \lambda_i^- )^{3/2}  < \infty$ and if $N^->0$, then $\lambda_1^- > 0$. 

\item the spectral measure $\mu$ of $J$ with decomposition $d\mu(x) = f(x)dx+d\mu_s(x)$ with respect to the Lebesgue measure satisfies
\begin{align*}
\int_{\tau^-}^{\tau^+} \frac{\sqrt{(\tau^+-x)(x-\tau^-)}}{x} \log(f(x)) dx >-\infty .
\end{align*}
\end{enumerate}
%
%
%
%
%
%
%
%
%
%
\end{cor}
\proof 
 It is enough to notice that $\mathcal F_L^- (0) = \infty$ and 
\[\mathcal F^\pm_L (\tau^\pm \pm h) = \frac{4}{3\tau^{3/4}(1 +\pm\sqrt \tau)^2 }h^{3/2} + o(h^{3/2}) \ \ \ (h \rightarrow 0^+)\]
and
\[G(1+h) = \frac{h^2}{2} + o(h) \ \ \ \ (h \rightarrow 0)\]
\begin{rem}
\label{remgemL}
Comparing with the Hermite gem, it would be most desirable to obtain a purely spectral criterion for when $J$ is a Hilbert-Schmidt operator relative to the Jacobi operator of the equilibrium measure, that is, in the Laguerre case
\begin{align} \label{HSlag}
\sum_{k=1}^\infty[(b_k-1-\tau)^2 + (a_k - \sqrt{\tau})^2] < \infty .
\end{align}
Unfortunately, Theorem \ref{sumrulel} will not yield such a criterion. (\ref{zl2}) implies (\ref{HSlag}) but the converse is not true. As an example, set $z_{2k-1}=\tau$, $z_{2k}=1$ for all $k\geq 1$, then (\ref{HSlag}) is clearly satisfied, but $\mathcal{I}_L(\mu)=\infty$ for $\mu$ the spectral measure of $J$ and $\tau \neq 1$. Actually, this system of coefficients correspond to the measure
\[\mu (dx) = (1-\tau) \delta_0 + \tau \MP_\tau (dx)\]
(see \cite{Saitoh} for the identification); the extra mass in $0$ gives a contribution $\mathcal F_L^- (0) = \infty$, the condition 2 is not fulfilled,
 although  the conditions 1 and 3 are fulfilled.
\end{rem}

\begin{cor} \label{semigemJ}
Assume the entries of the Jacobi matrix $J$ can be decomposed as in \eqref{zerl2} and let $\mu$ be the spectral measure of $J$ with pushforward $\tilde\mu$ under the mapping $x\mapsto \tfrac{1}{2} -\tfrac{1}{4}x$. Then, for any $\kappa_1,\kappa_2>0$, 
\begin{align}
\label{l2J}
\sum_{k=1}^\infty\left[\left( \alpha_{2k-1}+\frac{\kappa_1+\kappa_2}{2+\kappa_1+\kappa_2)}\right)^2 + \left( \alpha_{2k} - \frac{\kappa_1-\kappa_2}{2+\kappa_2+\kappa_2}\right)^2\right] < \infty
\end{align}
(that is, $\mathcal{I}_J(\tilde\mu)< \infty$) if and only if 
\begin{enumerate}
\item 
 $\tilde \mu \in \mathcal S_1(u_-, u_+)$
\item $\sum_{i=1}^{N^-} (\tfrac{1}{2}-\tfrac{1}{4}\lambda_i^-- u^+)^{3/2} + \sum_{i=1}^{N^+} (u^- -\tfrac{1}{2}+ \tfrac{1}{4}\lambda_i^+ )^{3/2}  < \infty$ and $\lambda_1^->-2$ if $N^->0$ and $\lambda_1^+<2$ if $N^+>0$.

\item when $d\tilde\mu(x) = f(x)dx+d\tilde\mu_s(x)$ is the decomposition of $\tilde\mu$ with respect to the Lebesgue measure, then
\begin{align*}
\int_{u^-}^{u^+} \frac{\sqrt{(u^+-x)(x-u^-)}}{x(1-x)} \log(f(x)) dx >-\infty .
\end{align*}
\end{enumerate}

\end{cor}

The proof is similar to the Laguerre case.

\begin{rem}
We can argue as in Remark \ref{remgemL}. In particular  (\ref{l2J}) implies
 \begin{align} \label{l2ab}
\sum_{k=1}^\infty[(\tilde b_k- \tilde b_k^{\KMK})^2 + (\tilde a_k - \tilde a_k^{\KMK})^2] < \infty\,,
\end{align}
but it is not equivalent.
\end{rem}

\section{Large deviations main theorem}
\label{sectLAR}

\subsection{The main result}
\label{sectLAR1}

Our large deviation result will hold for general eigenvalue distributions $\Pnv$ defined in \eqref{generaldensity}.
The corresponding spectral measure $\mun$ is then defined by \eqref{spectralmeasure}, where the weights $\w_1,\dots ,\w_n$ are Dir$_n(\beta')$ 
distributed and independent of the eigenvalues. We regard $\mun$ as a random element of $\mathcal{P}_1$, the set of all probability measures on $\R$, endowed with the weak topology and the corresponding $\sigma$-algebra. We need 
one more definition 
 in order to formulate the general result.

Recall that  $\mathcal{J}_V$ has been defined in  assumption (A3). We define, in the general case, the rate function for the extreme eigenvalues,
\begin{align}
\label{rate0}
\mathcal{F}_V^+(x) & = \begin{cases}
\mathcal{J}_V(x) - \inf_{\xi \in \R} \mathcal{J}_V(\xi) & \text{ if } \alpha^+\leq x \leq b^+, \\
\infty & \text{ otherwise, } 
\end{cases} \\
\mathcal{F}_V^-(x) & = \begin{cases}
\mathcal{J_V}(x) - \inf_{\xi \in \R} \mathcal{J}_V(\xi) & \text{ if } b^-\leq x \leq \alpha^-, \\
\infty & \text{ otherwise. } 
\end{cases}
\end{align}

\begin{thm} \label{MAIN}
Assume that the potential $V$ satisfies the assumptions (A1), (A2) and (A3). Then the sequence of spectral measures $\mun$ under $\Pnv\otimes \operatorname{Dir}_n(\beta')$ satisfies the LDP with speed $\beta'n$ and rate function
\begin{align*}
\mathcal{I}_V(\mu) = \mathcal{K}(\mu_V\!\ |\!\ \mu) + 
\sum_{n=1}^{N^+} {\mathcal F}_{V}^+ (\lambda_n^+)  +  \sum_{n=1}^{N^-} {\mathcal F}_{V}^- (\lambda_n^-)
\end{align*}
if $\mu \in \mathcal{S}_1(\alpha^-,\alpha^+)$ and $\mathcal{I}_V(\mu) = \infty$ otherwise. 
\end{thm}
Additionally, we have an alternative expression for $\mathcal F_V^\pm$, given by the following proposition. This result is
more or less classical. It may be found in \cite{Deiftuniform} (proof of Theorem 3.6) or in \cite{Albeverio} (Equation (1.13)).
\begin{prop}
\label{FandF}
If moreover, the conditions of analyticity (A5) and off-criticality (A4) are satisfied, then
\begin{align}
\mathcal F_{V}^+(x) & =   \int_{\alpha^+}^x S(t) \sqrt{(t-\alpha^-)(t-\alpha^+)}\!\ dt 
& \text{ if } x \geq \alpha^+ , \\ 
\mathcal F_{V}^-(x) & =  \int_x^{\alpha^-} S(t) \sqrt{(\alpha^- -t)(\alpha^+-t)}\!\ dt & \text{ if } x\leq \alpha^-.
\end{align}
\end{prop}

\subsection{From large deviations to sum rules}
\label{sectLAR2}

As described in the introduction, the sum rules of Section \ref{sectSU} are a consequence of two different proofs of a LDP, one leading to our main result Theorem \ref{MAIN}, giving the spectral side and another one yielding the sum side. Let us explain this in detail for the Hermite case. \\

For the Hermite case with probability measures on the whole real line, the correct randomization on the set of probability measures is the Hermite ensemble, defined by the eigenvalue density 
\begin{align*}
c_\beta \prod_{i<j}|\la_j-\la_i|^{\beta} \prod_{i=1}^n e^{-\frac{\beta' n}{2}\lambda_i^2}
\end{align*}
corresponding to the potential $V(x)= x^2/2$, and with weights following the Dirichlet distribution independent of the eigenvalues. Wigner's famous theorem states that the weak limit of the empirical eigenvalue distribution is then the semicircle law $\operatorname{SC}$. Indeed, here the potential $V$ satisfies all assumptions in Section \ref{sectLAR1} with $\mu_V=\operatorname{SC}$ and $S(x)=\tfrac{1}{2}$ on $[-2,2]$. Thus, by Theorem \ref{MAIN} the LDP for the measure $\mun$ holds. Further,  by Proposition \ref{FandF}, we may calculate the rate for the outliers as $\mathcal{F}_V^\pm= \mathcal{F}_H^\pm$. The rate function $\mathcal{I}_V$ is therefore precisely the right hand side in Theorem \ref{sumruleg}. \\
On the other hand, the recursion coefficients $(a_k)_k,(b_k)_k$ of the measure $\mun$ are independent with respectively gamma and normal distributions. Using this representation for the spectral measure, \cite{FGAR} proved that $\mun$ satisfies an LDP, again with speed $\beta' n$, and with rate function $\mathcal{I}_H$ the left hand side in Theorem \ref{sumruleg}. Since the rate function is unique, we must have $\mathcal{I}_V= \mathcal{I}_H$.\\

For the new sum rules, the arguments are similar. In the Laguerre case, the eigenvalue distribution of the spectral measure is
\begin{align*}
c_{\tau,\beta} \prod_{i<j} |\lambda_i-\lambda_j|^\beta \prod_{i=1}^n \lambda_i^{\beta' n (\tau^{-1}-1)} e^{- \beta' n \tau^{-1}\lambda_i} \mathbbm{1}_{\{ \lambda_i>0\} }
\end{align*}
with $\tau\in(0,1]$ and independent Dirichlet distributed weights. The potential of the Laguerre ensemble is $V(x) = \tau^{-1} x - (\tau^{-1} -1) \log x $ on $(0,\infty)$. As $n\to \infty$, the empirical eigenvalue distribution and the weighted spectral measure $\mun$ converge to the Marchenko-Pastur law $\operatorname{MP}_{\tau}$. Moreover, the assumptions of Theorem \ref{MAIN} are satisfied and we have an LDP with speed $\beta'n$ and rate function $\mathcal{I}_V$ the right hand side in the new sum rule, as $\mathcal{F}_V^\pm=\mathcal{F}_L^\pm$. As for the Hermite ensemble, \cite{FGAR} proved an LDP for $\mun$ in the subset of probability measures on $[0,\infty)$ with speed $\beta' n$ and rate function $\mathcal{I}_L(\mu) $ (note that \cite{FGAR} consider the speed $\beta' n \tau$). The uniqueness of the rate function implies $\mathcal{I}_V=\mathcal{I}_L$.\\

In the Jacobi case, the eigenvalue density is
\begin{align*}
c_{\kappa_1,\kappa_2,\beta} \cdot \prod_{i<j} |\lambda_i-\lambda_j|^\beta \prod_{i=1}^n \lambda_i^{\kappa_1\beta'n} (1-\lambda_i)^{\kappa_2\beta' n} \mathbbm{1}_{\{ 0<\lambda_i<1 \} }
\end{align*}
corresponding to the potential $V(x)= -\kappa_1\log(x)-\kappa_2\log(1-x)$ on $(0,1)$ for parameters $\kappa_1,\kappa_2\geq 0$. The equilibrium measure is then the Kesten-McKay distribution KMK$_{\kappa_1,\kappa_2}$. By Theorem \ref{MAIN}, $\mun$ satisfies the LDP with rate function $\mathcal{I}_V$, where additionally $\mathcal{F}_V^\pm=\mathcal{F}_J^\pm$. On the other hand, we know from the paper of \cite{FGAR}, that the LDP with rate function $\mathcal{I}_H$ holds. The combination of these two results yields Theorem \ref{sumrulej}.


\subsection{Conjecture for a general sum rule}
\label{conj}
\subsubsection{The probabilistic point of view}
\label{conjlaid}

We know from Section \ref{sectLAR1} that under  some assumptions on the potential  $V$, the random spectral measure sequence $(\mun)_n$ satisfies the LDP with rate function $\mathcal I_V$. Besides,  owing to Theorem \ref{KRV}, we can hope  to compute the rate function of the encoding by Jacobi coefficients directly from  the expression of the density. 
For a semi-infinite Jacobi matrix $ T=T((a_k)_k,(b_k)_k)$ with upper left $n\times n$ block $T\sn$ set 
\[\mathcal H (T\sn) = \tr\  V(T\sn) - 2\sum_{k=1}^{n-1} \log a_k \]
It would give for the rate function of the LDP at the speed $n\beta'$ 
\[T
\longmapsto 
\lim_{n\to \infty} \left[\mathcal H (T\sn)  - \inf_S \mathcal H(S\sn)\right]\]

 So, we conjecture the following identity, as soon as $V$ is a polynomial with even degree and positive leading coefficient :
\begin{eqnarray}
\label{general}
\mathcal K(\mu_V\!\ |\!\ \mu) + \sum_{k=1}^{N^+} \mathcal F_V(\lambda_k^+)+ \sum_{k=1}^{N^{-}} \mathcal F_V (\lambda_k^-)=\lim_{n\to \infty} \left[\mathcal H (J\sn_\mu)  - \inf_S \mathcal H(S\sn)\right]
\end{eqnarray}
Let us show that this is in agreement with the sum rules proven in this paper.
In the Hermite case when $V(x) = x^2/2$ , we get   
\[
\tr \ V(T\sn) = \frac{1}{2} \sum_{k=1}^n b_k^2 +  \sum_{k=1}^{n-1} a_k^2
\]
and then
\[
\mathcal H(T\sn) = \frac{1}{2}\sum_{k=1}^n b_k^2 + \sum_{k=1}^{n-1} a_k^2 - \log a_k^2\]
Now, $\inf_S \mathcal H(S\sn)$ is achieved for $b_k(S) \equiv 0, a_k(S) \equiv 1$, so that  
\[\lim_{n\to \infty} \left[\mathcal H (J\sn_\mu)  - \inf_S \mathcal H(S\sn)\right] =\sum_{k=1}^\infty  \frac{1}{2}b_k^2 + G(a_k^2) ,
\]
where $G(x) = x - 1 - \log x$ and this is exactly the rate function of (\ref{rateG1}). \\

In the Laguerre case when $V(x) = \tau^{-1} x - (\tau^{-1} -1) \log x$,
\[
\tr\ V(T\sn) = \tau^{-1} \tr\ T\sn - (\tau^{-1}-1) \log \det T\sn .
\]
But with the notation of Section \ref{susectLAG}, 
\[b_k=z_{2k-2}+z_{2k-1} \ \ \hbox{and} \ \  a_k = \sqrt{z_{2k-1}z_{2k}}\,\] so that $\tr\ T\sn = \sum_{k=1}^{2n-1} z_k$ and $T\sn = B\sn(B\sn)^*$, with 
\[
B\sn_{k,k}= \sqrt{z_{2k-1}},\quad B\sn_{k+1,k}= \sqrt{z_{2k}}
\]
and $B\sn_{i,j}=0$ for other entries. Then we have $\det T\sn = (\det B\sn)^2= \prod_{k=1}^{n} z_{2k-1}$ and
\begin{align*}
\mathcal H(T\sn) &= \tau^{-1} \sum_{k=1}^{2n-2} z_k -  (\tau^{-1}-1)\sum_{k=1}^n \log z_{2k-1} -  \sum_{k=1}^{2n-1} \log z_k \\
&= \tau^{-1}\sum_{k=1}^{n} (G(z_{2k-1}) +1) + \sum_{k=1}^{n-1} (G(\tau^{-1}z_{2k})-\log(\tau)+1) + \log z_{2n-1} .
\end{align*}
Of course, the infimum is achieved for $z_{2k-1}\equiv 1, z_{2k} \equiv \tau$, which gives exactly the expression of (\ref{rateL1}).\\

Finally, let us look at the Jacobi case, transformed onto the interval $[-2,2]$, with \[V(x) = -\kappa_1 \log (2-x) - \kappa_2 \log (2+x)\] and
\begin{align*}
\tr\ V(T\sn) =-\kappa_1 \log \det (2I\sn-T\sn) - \kappa_2 \log \det (2I\sn+T\sn) .
\end{align*}
In formula (5.2) of \cite{Killip1} we see
\[\Phi_{2n}(1) = \prod_{k=0}^{2n-1} (1-\alpha_k) = \prod_{k=1}^n (2-\lambda_j) = \det (2I\sn-T\sn)
\]
Similarly, from formula (5.3) of the same paper we get
\[
\Phi_{2n}(-1) = \prod_{k=0}^{2n-1} (1 + (-1)^{k} \alpha_k) = \prod_{k=1}^n (2+\lambda_j) = \det (2I\sn+T\sn) 
\]
such that we obtain
\begin{align*}
\tr \ V(T\sn) & = -\kappa_1 \log \prod_{k=0}^{n-1}(1-\alpha_{2k})
-\kappa_2 \log \prod_{k=0}^{n-1}(1+\alpha_{2k})\\
& \quad 
-\kappa_1 \log \prod_{k=0}^{n-1}(1-\alpha_{2k+1})
-\kappa_2 \log \prod_{k=0}^{n-1}(1-\alpha_{2k+1})
\end{align*}

Recall that $a_{k+1}^2= (1-\alpha_{2k-1})(1-(\alpha_{2k})^2)(1+\alpha_{2k+1})$, and then
\begin{align*}
\mathcal{H}(T\sn)&= \tr \ V(T\sn) - \sum_{k=0}^{n-1} \log(1-\alpha_k^2) + A_n\\
&= -\sum_{k=0}^{n-1} (1+\kappa_1) \log(1-\alpha_{2k}) + (1+\kappa_2)\log(1+\alpha_{2k}) \\
& \quad - \sum_{k=0}^{n-1} (1+\kappa_1+\kappa_2)\log(1-\alpha_{2k+1})+ \log(1+\alpha_{2k+1}) 
 \quad +A_n
\end{align*}
with $A_n=\log(1-\alpha_{2n-3}) +\log(1-\alpha_{2n-2}^2)+\log(1-\alpha_{2n-1}^2)$. Again, in the limit $n\to \infty$ this leads exactly to the rate function \eqref{rateJ1}. 

\subsubsection{Mathematical analysis point of view}
\label{analysis}
In \cite{Nazarov}, a sum rule is given when the reference measure may be written as
\begin{equation}
\label{naz1}\sigma(dx) = A(x) \operatorname{SC}(dx)\end{equation}
with $A$ a nonnegative polynomial (Theorem 1.5 therein). Under appropriate conditions, the sum rule is 
\begin{equation}
\label{Nazarov}
\mathcal K(\sigma\!\ |\!\ \mu) + \sum_{k=1}^{N^+} \mathcal F(\lambda_k^+)+ \sum_{k=1}^{N^{-}} \mathcal F (\lambda_k^-)=\lim_{n \to \infty} \left(- 2 \sum_1^{n-1} \log a_k + \tr \left(\Phi(T\sn) - \Phi(T_0\sn)\right)\right)\,
\end{equation}  
where 
\[\mathcal F (x) =  \begin{cases} \displaystyle \int_2^x A(t) \sqrt{t^2-4}\!\ dt \;\;\;\;& \mbox{if}\ x \geq 2\,, \\
\displaystyle \int_x^{-2} A(t) \sqrt{t^2-4}\!\ dt  &\mbox{if}\ x \leq -2\,.
\end{cases}\]
Actually, these conditions warrant the existence of this limit.  Set
\begin{equation}
\label{aux}\Phi'(z) = zA(z) - \frac{1}{\pi} \int \frac{A(x) - A(z)}{x-z} \operatorname{SC}(dx)\end{equation}
which leads for $\  z \notin [-2,2]$  to
\[\Phi'(z) = \mathcal F'(z) - \int\frac{\sigma(dx)}{x-z}
\,.\]
The function $\Phi$ was defined in  \cite{Nazarov} as 
an auxiliary function.
In view of our Proposition \ref{FandF} (up to an affine change of variables), it  appears then that the triple $(\sigma, \mathcal F, \Phi)$ is actually identical to the triple $(\mu_V, \mathcal F_V, V)$.

Nazarov et al.  claimed that the scope of generality was not clear to them, so that they use the polynomial nature of $A$. 
Actually, our classical ensembles are, again up to an affine change of variables, of the form (\ref{naz1}) with $A$ analytic in a neighboorhood of $[-2,2]$, namely with $1/A$ polynomial of degree at most 2.

To end this section let us give some concluding remarks.

\begin{enumerate}
\item The sum rule given by  Nazarov et al. should be true  beyond  the polynomial case for more general functions $A$. 
\item The function $\Phi$ is nothing else than a potential and $\mathcal F$  the corresponding effective potential.
\item When $A$ is a polynomial the underlying potential is also a polynomial.
\item In this latter case, Section \ref{conjlaid} gives  a probabilistic interpretation of (\ref{Nazarov}) and a draft for a probabilistic proof. 
\end{enumerate}
We also refer to \cite{KuKu} for other extensions.

\section{Proofs}
\label{sectproof}

This section is devoted to the proof of Theorem \ref{MAIN}. The main idea is to apply the projective limit method to reduce the spectral measure to a measure with only a fixed number of eigenvalues outside the limit support $[\alpha^-,\alpha^+]$. 
For this we need to consider a topology on $\Sr$ different from the weak topology. Recall that measures $\mu \in \Sr$ are written as
\begin{align}\label{muinS}
\mu = \mu_{|I} +  \sum_{i=1}^{N^+} \gamma_i^+ \delta_{\lambda_i^+} + \sum_{i=1}^{N^-} \gamma_i^- \delta_{\lambda_i^-}
\end{align}
and we associate the measure \eqref{muinS} with 
\begin{align}\label{muinS2}
\big( \mu_{|I}, (\lambda_i^+)_{i \geq 1}, (\lambda^-_i)_{i\geq 1},(\gamma^+_i)_{i\geq 1},(\gamma^-_i)_{i\geq 1}\big)
\end{align}
 with $\lambda_i^+=\ap$ and $\gamma_i^+=0$ if $i>N^+$ and $\lambda^-_i=\am$ and $\gamma_i^-=0$ if $i>N^-$. 
The topology on $\Sr$ is then defined by the vector \eqref{muinS2}: we say that 
 $\mu_n$ converges to $\mu$ in $\Sr$ if:
\begin{align} \label{strangetop}
\begin{split} \mu_{n|I} \xrightarrow[n \rightarrow \infty ]{} \mu_{|I} & \text{ weakly and 
for every}\  i \geq 1 \\
\big( \lambda^+_i (\mu_n),\lambda^-_i(\mu_n),\gamma^+_i(\mu_n),\gamma^-_i&(\mu_n)\big) \xrightarrow[n \rightarrow \infty ]{}\big( \lambda^+_i (\mu),\lambda^-_i(\mu),\gamma^+_i(\mu),\gamma^-_i(\mu)\big) 
\end{split}
\end{align}
We will show in Section \ref{susectPROJL} that on the smaller set $\Sr_1=\{\mu \in \Sr |\, \mu(\R)=1\}$, 
this convergence implies weak convergence, but we remark that $\mu_n \to \mu$ weakly does not imply convergence in our topology. For example, the merging of two atoms outside of $I$ is no continuous operation, while it is continuous in the weak topology. The $\sigma$-algebra on $\Sr$ is then the corresponding Borel-algebra. \\
On $\Sr$ we define a family of projections $(\pi_j)_{j\in\mathbb{N}}$, where for a measure $\mu$ as in \eqref{muinS}, 
\begin{align}
\label{defproj}
\pi_j(\mu) = \mu_{|I} + \sum_{i=1}^{N^+\wedge j} \gamma_i^+ \delta_{\lambda_i^+} + \sum_{i=1}^{N^-\wedge j} \gamma_i^- \delta_{\lambda_i^-}\,,
\end{align} 
that is, we keep $\mu$ on $I$ but delete all but up to $j$ point masses left of $\am$ and right of $\ap$. Note that the projections are continuous in our topology, but they are not in the weak topology.

\label{sectSU2}
\subsection{LDP for a finite collection of extreme eigenvalues}

The study of LDP for (one) extreme eigenvalue of random matrices began in \cite{aging} and in \cite{Albeverio}. For detailed comments on the assumptions see Section \ref{comm}.

\subsubsection{Notation}
Under the probability measures considered, there are almost surely no ties among eigenvalues, so that 
we may reorder  $\lambda = (\lambda_1, \dots, \lambda_n)$ as $\hat\lambda= (\lambda_1^n, \dots, \lambda_n^n)$ such that $\lambda_1^n > \lambda_2^n > \dots >\lambda_n^n$. Let $\lambda^+_i=\lambda^n_i$ and $\lambda^-_i=\lambda^n_{n-i+1}$ and 
for $j$ a fixed integer and
$n > 2j$ 
\[ \lambda^+(j) = ( \lambda_1^+, \dots, \lambda_j^+)\ , \ \lambda^-(j) = (\lambda_1^-, \dots, \lambda^-_{j})\,.\]

For the sake of simplicity, we denote by $\mathbb R^{\uparrow j}$ (resp. $\mathbb R^{\downarrow j}$) the subset of $\mathbb R^j$ of all vectors with non decreasing (resp. non increasing) components.
\subsubsection{Main result}
\begin{thm}
\label{LDPjextreme}
Let $j$ and $\ell$ be  fixed integers. Assume that $V$ is continuous and satisfies (A1), (A2) and the control condition (A3).
\begin{enumerate}
\item If $b^-<\alpha^-$ and $\alpha^+<b^+$, then
the law of $(\lapj, \lambda^-(\ell))$ under $\Pnv$ satisfies a LDP in $\R^{j+ \ell}$ with speed $\beta'n$ and rate function
\begin{align*}
\mathcal{I}_{\lambda^\pm}(x^+, x^-)
 := 
\begin{cases} &  
\sum_{k=1}^j 
\Fr_{V}^+(x_k^+)  + \sum_{k=1}^\ell 
\Fr_{V}^-(x^-_k)\
%
\;\;\;\;\mbox{if}\ 
(x_1^+, \dots , x_j^+)\in \mathbb R^{\downarrow j} \ \hbox{and}\  (x_1^-, \dots , x_\ell^-)\in \mathbb R^{\uparrow \ell}\\
      &  \infty\;\;\mbox{ otherwise.}
\end{cases}
%
\end{align*}
\item If $b^-=\alpha^-$, but $\alpha^+<b^+$, the law of $\lapj$ satisfies the LDP with speed $\beta'n$ and rate function
\begin{align*}
\mathcal{I}_{\lambda^+}(x^+)= 
\mathcal{I}_{\lambda^\pm}(x^+, \alpha^-)=
\begin{cases} &  
\sum_{k=1}^j 
\Fr_{V}^+(x_k^+)
\;\;\;\;\mbox{if}\ 
(x_1^+, \dots , x_j^+)\in \mathbb R^{\downarrow j}\\
      &  \infty\;\;\mbox{ otherwise.}
\end{cases}
\end{align*}
\item If $b^-<\alpha^-$, but $\alpha^+=b^+$, the law of $\lambda^-(\ell)$ satisfies the LDP with speed $\beta'n$ and rate function
\begin{align*}
\mathcal{I}_{\lambda^-}( x^-) = \mathcal{I}_{\lambda^\pm}(\alpha^+, x^-) =
\begin{cases} &
\sum_{k=1}^\ell 
\Fr_{V}^-(x^-_k)\
\;\;\;\;\mbox{if}\   (x_1^-, \dots , x_\ell^-)\in \mathbb R^{\uparrow \ell}\\
     &  \infty\;\;\mbox{ otherwise.}
\end{cases}
\end{align*}
\end{enumerate}
\end{thm}

The same statement is Theorem 2.10 in \cite{benaych2012large}, but with an extra technical assumption that is not easy to check.
Besides, after the above mentioned publications, recently  
\cite{Borot} and \cite{borot-multi} provided other sketch of proofs for the case where $\ell=0, j=1$ and without this assumption. 
For the sake of completeness, we give a  proof for the general case.  As this proof  is technical, we postpone it to the  Appendix (Section \ref{Appendix1}).

\label{susectEXT}
\subsection{Joint LDP for the restricted measure and a finite collection of extreme eigenvalues}
\label{susectJOIN}

Recall that $\left((\lambda_1, \dots, \lambda_n), (\w_1, \dots , \w_n)\right)$ is distributed according to $\mathbb Q_n^V = \Pnv \otimes \hbox{Dir}_n(\beta')$. 
The two sources of randomness are not at the same scale. On the one hand,  the eigenvalues are ruled by a LDP at speed $n^2$. On the other hand,  the weights are ruled by a factor $n$. Hence,  it is natural to consider the eigenvalues as quasi-deterministic, and to begin by a conditioning upon these variables. As in a previous paper (\cite{gamboacanonical}), it is then convenient to 
decouple the weights by introducing independent random variables. We know that
\begin{equation}
\label{betagamma}(\w_1, \dots, \w_n) \el \left(\frac{\gamma_1}{\gamma_1+ \dots + \gamma_n}, \dots , \frac{\gamma_n}{\gamma_1+ \dots + \gamma_n}\right)\end{equation}
where $\el$ means equality in distribution, and $\gamma_1, \dots, \gamma_n$ are independent variables 
with distribution Gamma$(\beta', (\beta'n)^{-1})$ and mean $n^{-1}$. We enlarge the probability space to define such variables $\gamma_i$'s and denote by $\widetilde {\mathbb Q_n^V}$ the corresponding probability measure. With this notation, we can rewrite the spectral measure $\mu\sn$ as
\begin{equation}
\label{normunnorm}\mu\sn = \frac{\widetilde \mu\sn}{\widetilde \mu\sn(\mathbb R)}\end{equation}
where 
\begin{equation}
\label{decouple}\widetilde \mu\sn := \sum_{k=1}^n \gamma_k \delta_{\lambda_k}\end{equation}
is a random measure with independent masses $\gamma_1, \dots, \gamma_n$ and $\widetilde \mu\sn(\mathbb R)$ is its total mass $\sum_{k=1}^n \gamma_k$.  

We denote by $\munI$ the restriction of $\mun$ to the interval $I$.
 Similarly, for $I(j)=I\setminus \{\lap_1,\lambda_1^-,\dots ,\lap_j,\lambda_j^- \}$, $\munIj$ is the restriction of  $\mun$ to $I(j)$ and we use the analogous notation for the restrictions of the empirical measure $\muun$.
 Notice that we choose $j=\ell$ for the sake of simplicity.  
 The aim of this subsection is to prove the following joint LDP for the restricted spectral measure and a collection of largest and/or smallest eigenvalues.

\medskip

\begin{thm} \label{jointLDP}
$\;$
\begin{enumerate}
\item If $\bm<\am<\ap<\bp$, then for any fixed $j\in\N$ and under $\widetilde{\mathbb Q_n^V}$, the sequence of random objects
$\big( \mutnIj, \lapj, \lamj \big)$
satisfies the joint LDP with speed $\beta'n$ and rate function
\begin{align*}
\mathcal{I}(\mu,x^+,x^-) =  \mathcal{K}(\mu_V\!\ | \!\  \mu) +\mu(I) - 1+ \mathcal{I}_{\lambda^\pm}(x^+, x^-)
\end{align*}
\item If $\bm=\am$, but $\ap<\bp$ (or $\bp=\ap$, but $\am>\bm$), then, with the same notation as in the previous section, 
$
\big( \munIj, \lapj \big) (\text{or } \big( \munIj, \lamj  \big) \text{ respectively,} )
$
satisfies the LDP with speed $\beta'n$ and rate function 
\begin{align*}
\mathcal{I}^+(\mu,x^+) = \mathcal{I}(\mu,x^+,\am) \quad (\text{or } \ \mathcal{I}^-(\mu,x^-) = \mathcal{I}(\mu,\ap,x^-) \text{ respectively})\,.
\end{align*}
\end{enumerate}
\end{thm}

\medskip

\proof
We only prove here the first point of the theorem. The second claim can be shown in the same way.
We first show a joint LDP when the eigenvalues are truncated. For $M>\max\{|\ap|,|\am|\}$, let $\lapjM$ (resp.$\lamjM$) be the collection of truncated eigenvalues 
\begin{align*}
\lambda_{M,i}^+ = \min\{ \lambda_i^+ , M\}\ \ \ (\text{resp.}\  \lambda_{M,i}^- = \max\{ \lambda_i^- , -M\})\,,
\end{align*}
for $i=1,\dots,j$. 
To further simplify notation, let $\lambda_M^\pm(j)=(\lambda_{M,
1}^+,\dots,\lambda_{M,j}^+,\lambda_{M,1}^-,\dots ,\lambda_{M,j}^+)$. \\

\textbf{Exponential tightness:}\\
In a first step, we will obtain the joint LDP for $( \mutnIj, \lambda_M^\pm(j))$ by applying Theorem 1.1 of \cite{baldi1988large}. For this, we need to check that this sequence is exponentially tight. For $M$ as above, define the set
\begin{align*}
K_M = \left\{ (\mu,\lambda) \in \mathcal{S}\times \mathbb{R}^{2j} |\, \mu(I) \leq M, \mu(I^c)=0, \lambda \in [-M,M]^{2j} \right\} .
\end{align*}
Indeed, $K_M$ is a compact set in the topology \eqref{strangetop} and
\begin{align*}
\mathbb{P} ( ( \mutnIj, \lambda_M^\pm(j)) \notin K_M ) = \mathbb{P} (\mutnIj(I) > M) 
 \leq \mathbb{P} \left( \sum_{k=1}^n \gamma_k >M \right) .
\end{align*}
The sum in the last probability is Gamma$(\beta'n, (\beta'n)^{-1})$ distributed. By the LDP for the Gamma-distribution with rate $G$,  
\begin{align*}
\mathbb{P} \left( \sum_{k=1}^n \gamma_k >M \right)  \leq e^{-\beta'n G(M)}
\end{align*}
Therefore,
\begin{align*}
\limsup_{n\to \infty} \frac{1}{\beta' n} \log \mathbb{P} ( ( \mutnIj, \lambda_M^\pm(j)) \notin K_M ) \leq - G(M) ,
\end{align*}
which can be chosen to be arbitrarily small, i.e., the sequence $( \mutnIj, \lambda_M^\pm(j))$ is exponentially tight.\\

\textbf{Joint LDP for measure and truncated eigenvalues:}\\
Let $f$ be a continuous function from $\mathbb R$ to $\mathbb R$ such that $\log(1-f)$ is bounded.  For $s^\pm \in\R^{2j}$, we calculate the joint moment generating function
\begin{align*}
&\mathcal{G}_n(f,s^\pm) = 
\E\left[ \exp\left\{ n \beta' \left(\mutnIj(f) + \langle s^\pm,\lambda_M^\pm(j)\rangle \right) \right\} \right] 
\end{align*}
under $\widetilde{\mathbb Q_n^V}$.
First recall that  $\gamma_i$ is $\operatorname{Gamma}(\beta',(\beta'n)^{-1})$ distributed, so that 
\begin{align}
\label{cgfgamma}
\frac{1}{\beta'}\log \E e^{\beta'n\gamma_i t} =  L(t)
\end{align}
(see (\ref{Laplace2}))  and then, 
 integrating with respect to the $\gamma_i$'s  we get
\begin{align*}
\mathcal{G}_n(f,s^\pm) &=   \mathbb E \left[\exp\left(n\beta' \langle s^\pm , \lambda_M^\pm (j)\rangle\right) \prod_{i \in I(j)} \mathbb E\left[e^{n\beta'\gamma_1 f(\lambda_i)}| \lambda_1, \dots, \lambda_n\right]\right]\\
&= 
   \E\left[ \exp\left\{ n\beta' \left(\muunIj(L \circ f) + \langle s^\pm,\lambda_M^\pm(j)\rangle  \right)\right\} \right] \,,
\end{align*}
This expectation only involves $\Pnv$.
Set
\[D_n(s^\pm):=  \E \left[\exp\left\{ n\beta'\langle s^\pm,\lambda_M^\pm(j)\rangle \right\} \right]\,.\]
By Theorem \ref{LDPjextreme} we have a LDP for the extremal eigenvalues $\lambda^\pm(j)$ of the spectral measure with rate function $\mathcal{I}_{\lambda^\pm}$. By the contraction principle (see \cite{demboz98} p.126), the truncated eigenvalues satisfy the LDP with rate function
\begin{align*}
\mathcal{I}_{M,\lambda^\pm}(x^\pm) = \begin{cases} \mathcal{I}_{\lambda^\pm}(x^+,x^-) & \text{ if } x^\pm = (x^+,x^-) \in [-M,M]^{2j}, \\ \infty & \text{ otherwise.} \end{cases}
\end{align*}
Since the truncated eigenvalues are bounded, 
 Varadhan's Integral Lemma (\cite{demboz98} p. 137) implies 
\begin{align} \label{jointLDPproof1}
\lim_{n\to \infty} \frac{1}{\beta'n}\log D_n (s^\pm)
 = \mathcal{I}_{M,\lambda^\pm}^*(s^\pm) ,
\end{align}
where
\begin{align*}
\mathcal{I}_{M,\lambda^\pm}^* (s^\pm)= \sup_{x^\pm \in \R^{2j}} \left\{ \langle s^\pm,x^\pm\rangle   - \mathcal{I}_{M,\lambda^\pm}(x^\pm) \right\}
\end{align*}
is the convex dual of $\mathcal{I}_{M,\lambda^\pm}$. To control $ n \muunIj(L \circ f)$, let for $\eta>0$
\begin{align*}
 A(\eta)=\left\{d (\muunIj , \mu_V)< \eta \right\} ,
\end{align*}
with a metric $d$ inducing weak convergence. Since $\muunIj$ and $\muunI$ differ only by at most $2j$ support points, their total variation distance is bounded by $2j/{n}$. For $n$ large enough this implies
\begin{align*}
\left\{d( \muunI,  \mu_V)< \eta/2 \right\} \subset A(\eta) .
\end{align*}
Now,
\[\mathbb P(A(\eta)^c) \leq \mathbb P(d(\muunI, \mu_V)\geq \eta/2)\leq \mathbb P(d(\muun,\mu_V)\geq \eta/2)\]
and then, since 
 $\muun$ satisfies an LDP with speed $n^2$ and a rate which is good with unique minimizer $\mu_V$ (Theorem \ref{LDPmuu}) we have for $n$ large enough
\begin{align} \label{jointLDPproof2}
\mathbb P(A(\eta)^c) \leq 
e^{-n^2\delta} 
\end{align}
with a $\delta=\delta(\eta)>0$. Writing $\mathcal{G}_n(f,s^\pm)=\mathcal{G}_{n,A}(f,s^\pm)+\mathcal{G}_{n,A^c}(f,s^\pm)$ with 
\begin{align*}
\mathcal{G}_{n,A}(f,s^\pm) = E\left[ \exp\left\{ n\beta'\left( \muunIj(L \circ f) + \langle s^\pm,\lambda_M^\pm(j)\rangle \right) \right\} \indi_{A(\eta)}\right]
\end{align*}
we can bound
\begin{align}
\label{updownA} 
 C_n(s^\pm) \exp\left\{ n\beta'\left( \mu_V(L \circ f)-\eta \right)\right\} \leq \mathcal{G}_{n,A}(f,s^\pm)\leq  C_n(s^\pm)  \exp\left\{ n\beta'\left( \mu_V(L \circ f)+\eta \right)\right\}
\end{align}
where
\[C_n (s^\pm):=  \E \left[\exp\left\{ n\beta'\left(\langle s^\pm,\lambda_M^\pm(j)\rangle \right) \right\} \indi_{A(\eta)}\right] \leq D_n(s^\pm)\,,\]
and then from (\ref{jointLDPproof1})
\begin{align}
\label{upA}
\limsup _{n\to \infty} \frac{1}{\beta'n} \log \mathcal{G}_{n,A}(f,s^\pm) \leq  \mu_V(L \circ f)+\eta +\mathcal{I}_{M,\lambda^\pm}^*(s^\pm)\,.
\end{align}
For the complimentary event, we have the upper bound
\begin{align*}
\mathcal{G}_{n,A^c}(f,s^\pm) \leq  (D_n(s^\pm) -C_n(s^\pm)) \exp\left\{ n\beta' ||L \circ f||_\infty \right\}\,.
\end{align*}
By the Cauchy-Schwarz inequality and (\ref{jointLDPproof2}) we get
\[ (D_n(s^\pm) -C_n(s^\pm)) \leq D_n (2s^\pm) e^{-n^2\delta}\,,\]
and then, using (\ref{jointLDPproof1}) 
   we get
\begin{align}
\label{upD-C}
\limsup_{n\to \infty} \frac{1}{n\beta'} \log(D_n(s^\pm) -C_n(s^\pm))  = -\infty
\end{align}
which eventually leads to
\begin{align}
\label{upAc}
\limsup_{n\to \infty} \frac{1}{n\beta'} \log \mathcal{G}_{n,A^c}(f,s^\pm) = -\infty .
\end{align}
Combining \eqref{upA} and \eqref{upAc}, we get
\begin{align}
\label{ls}
\limsup_{n\to \infty} \frac{1}{n\beta'} \log \mathcal{G}_{n}(f,s^\pm) \leq \mu_V(L \circ f)+\eta + \mathcal{I}_{M,\lambda^\pm}^*(s^\pm) .
\end{align}
For the lower bound, we have
\begin{align*}
 & \liminf_{n\to \infty} \frac{1}{n\beta'} \log \mathcal{G}_{n}(f,s^\pm) \geq \liminf_{n\to \infty} \frac{1}{n\beta'} \log \mathcal{G}_{n,A}(f,s^\pm)\\
& \geq  \mu_V(L \circ f)- \eta  + \liminf_{n\to \infty} \frac{1}{n\beta'} \log C_n(s^\pm)\,,
\end{align*}
and from \eqref{jointLDPproof1} and \eqref{upD-C} 
\[\lim_{n\to \infty}  \frac{1}{n\beta'} \log  C_n(s^\pm) =  \mathcal{I}_{M,\lambda^\pm}^*(s^\pm)\,,\]
so that
\begin{align}
\label{li}
\liminf_{n\to \infty} \frac{1}{n\beta'} \log \mathcal{G}_{n}(f,s^\pm) \geq \mu_V(L \circ f)- \eta +  \mathcal{I}_{M,\lambda^\pm}^*(s^\pm)\,.
\end{align}
From (\ref{ls}) and (\ref{li}) and 	since $\eta>0$ was arbitrary, we can conclude:
\begin{align*}
\lim_{n\to \infty} \frac{1}{\beta'n} \log \mathcal{G}_{n}(f,s^\pm) = \mu_V(L \circ f) + \mathcal{I}_{M,\lambda^\pm}^*(s^\pm) =: \mathcal{G} (f,s^\pm) .
\end{align*}
The convex dual of $\mathcal{G}$ is
\begin{align*}
\mathcal{G}^* &(\mu,\lambda^\pm) = \sup_{f \in C_b(I),s^\pm \in \R^{2j}} \left( \int f d\mu +\langle \lambda^\pm,s^\pm\rangle  -\mathcal{G}(f,s^\pm)\right) \\
& = \sup_{f \in C_b(I)} \left( \int f d\mu - \mu_V(L \circ f) \right) + \sup_{s^\pm \in \R^{2j}} \left( \langle \lambda^\pm,s^\pm\rangle -\mathcal{I}_{M,\lambda^\pm}^*(s^\pm) \right)\\
& = \Lambda^*(\mu) + \mathcal{I}_{M,\lambda^\pm}(\lambda^\pm) ,
\end{align*}
where 
 $\Lambda^*$ 
is the convex dual of $\mu_V(L\circ \cdot)$. The LDP follows now from Theorem 1.1 of \cite{baldi1988large} with rate function given by $\mathcal{G}^*$, provided that $\mathcal{G}^*$ is strictly convex on a set of points that is dense in the set of all points where 
$\mathcal{G}^*$ is finite. \\

\textbf{
Identification of $\Lambda^*$:}\\
By Theorem 5 of  \cite{Rocky1}, we can write $\Lambda^*$ as
\begin{align*}
\Lambda^* (\mu) = \mu_V(L^* \circ h_\mu) +r(1) \mu_s(\R) ,
\end{align*}
where $L^*$ is the convex dual of $L$ and $r$ its recession function and $d\mu=h_\mu\cdot d\mu_V + d\mu_s$ is the Lebesgue-decomposition of $\mu$ with respect to $\mu_V$.
The expression  of $L^*$ is given in (\ref{Laplace2}).
 The recession function is 
\begin{align*}
r(x) = \sup \{ xy\ | \ L(y)<\infty\} = x .
\end{align*}
for nonnegative $x$. We obtain 
\begin{align*}
\Lambda^*(\mu) &=  -\int \log \left( h_\mu\right) d\mu_V  - 1 + \int h_\mu \!\ d\mu_V + \mu_s(\R)\\
& =  - \int \log \left( h_\mu\right) d\mu_V  - 1 + \mu(I) \\
& = \mathcal{K}(\mu_V | \mu) -1 + \mu(I)\,.
\end{align*}
Now $\Lambda^*$ is strictly convex at $\mu$ if there exists a $f\in C_b(I)$, called an exposing hyperplane, such that
\begin{align} \label{jointLDPproof3}
\Lambda^*(\mu)-\int fd\mu< \Lambda^*(\nu)-\int fd\nu
\end{align}
for any $\nu \neq \mu$. Suppose $d\mu =h_\mu\cdot d\mu_V$ is absolutely continuous with respect to $\mu_V$ with a density $h_\mu$ positive on $I$ and choose $f = 1- h_\mu^{-1}$. Then \eqref{jointLDPproof3} is equivalent to 
\begin{align*}
\int \log(h_\mu/h_\nu)d\mu_V > \int d\mu_V - \int h_\mu^{-1}d\nu,
\end{align*}
The last inequality follows from $\log(x)>1-x^{-1}$ , ($x>0, x\neq 1$). Indeed, 
\begin{align*}
\int \log(h_\mu/h_\nu)d\mu_V \geq \int (1-h_\mu^{-1}\cdot h_\nu) d\mu_V\geq \int d\mu_V - \int h_\mu^{-1} d\nu , 
\end{align*}
where the first inequality is strict unless $\mu_V$-almost everywhere $h=h_\nu$ and the second one is strict unless $\nu_s=0$. So if $\nu\neq \mu$ at least one inequality is strict, i.e. $\Lambda^*$ is strictly convex at all points $d\mu=h\cdot d\mu_V$, which are dense in the set of nonnegative measures on $I$. Consequently, $( \mutnIj, \lambda_M^\pm(j) )  $ satisfies an LDP with speed $\beta'n$ and rate function
\begin{align*}
\mathcal{I}(\mu,x^\pm) = \mathcal{K}(\mu_V\!\ |\!\ \mu) + \mu(I) - 1 +  \mathcal{I}_{M,\lambda^\pm}(x^\pm) .
\end{align*}

\textbf{Extending the LDP to untruncated eigenvalues:}\\
From the LDP for $(\lapj, \lamj )$ the exponential tightness of the (unrestricted) extremal eigenvalues holds (see Section \ref{expotight}). 
This implies
\begin{align*}
\lim_{M\to \infty}  \limsup_{n\to \infty} \frac{1}{n} \log P\big( \lambda_M^\pm(j) \neq \lambda^\pm(j) \big) = -\infty ,
\end{align*}
so that as $M\to \infty$, the truncated eigenvalues are exponentially good approximation of the unrestricted ones. In fact, $(\mutnIj,\lambda_M^\pm(j) )$ are exponentially good approximations of $(\mutnIj,\lambda^\pm(j) ) $. Since the rate function of the untruncated eigenvalues can be recovered as the pointwise limit
\begin{align*}
\mathcal{I}_{\lambda^\pm}(\lambda^\pm) = \lim_{M\to \infty} \mathcal{I}_{M,\lambda^\pm}(\lambda^\pm) ,
\end{align*} 
we get from Theorem 4.2.16 in \cite{demboz98} that 
$(\mutnIj,\lambda^\pm(j) ) $ satisfies the LDP with speed $\beta'n$ and rate function 
\begin{align*}
\mathcal{I}(\mu,x^\pm) = \mathcal{K}(\mu_V\!\ |\!\ \mu) + \mu(I) - 1 + \mathcal{I}_{\lambda^\pm}(x^\pm) = \mathcal{K}(\mu_V\!\ |\!\ \mu) + \mu(I) - 1 +  \sum_{i=1}^ {j} \Fr^+(x^+_i) + \Fr^-(x^-_i)\,,
\end{align*}
which ends the proof of Theorem \ref{jointLDP}.
\QED
\subsection{LDP for the projected measure}
\label{susectPROJ}
Recall the definition of the projections $\pi_j$ in (\ref{defproj}).
\begin{thm}\label{projLDP}
For any fixed $j$, the sequence of projected spectral measures $\pi_j(\mutn)$ as elements of $\mathcal{S}$ with topology \eqref{strangetop} satisfies under $\widetilde{\mathbb Q_n^V}$ the LDP with speed $\beta'n$ and rate function
\begin{align*}
\tilde{\mathcal{I}}_j(\tilde{\mu}) = \mathcal{K}(\mu_V \ | \ \tilde{\mu}) +\tilde{\mu}(I)-1+ \sum_{i=1}^{N^+\wedge j} \left(\Fr^+_V(\lambda^+_i)+\gamma_i^+\right) + \sum_{i=1}^{N^-\wedge j}\left(\Fr^-_V(\lambda^-_i)+\gamma_i^-\right).
\end{align*}
\end{thm}

\medskip

\proof
This result is a direct consequence of  Theorem \ref{jointLDP} and the contraction principle. 
Again, suppose that the first case of Theorem \ref{jointLDP} holds, otherwise $N^+$ or $N^-$ is 0 and we may omit the largest and/or smallest eigenvalues and weights. 
By the independence of the weights $\gamma_k$ and independence of weights and eigenvalues, 
the collection of weights $(\gapj,\gamj)$ is independent of $( \mutnIj, \lapj, \lamj )$. Recall that $(\gapj,\gamj)$ contains a collection of independent $\operatorname{Gamma}(\beta',(\beta'n)^{-1})$ distributed random variables and then satisfies an LDP with speed $\beta'n$ and rate function 
\begin{align*}
\mathcal{I}_\gamma (y^+,y^-) = \sum_{i=1}^j (y^+_i+y^-_i)
\end{align*}
for $y^+_i,y^-_i\geq 0$ and $\mathcal{I}_\gamma (y^+,y^-)=\infty$ otherwise (see (\ref{LDPdeg})). Thus, 
\begin{align}
\big( \mutnIj, \lapj, \lamj, \gapj,\gamj \big)
\end{align}
satisfies the LDP with speed $\beta'n$ and rate function
\begin{align*}
\mathcal{I}(\mu,x^+,x^-,y^+,y^-) = \mathcal{K}(\mu_V\!\ |\!\ \mu) +\mu(I) -1+ \sum_{i=1}^j \Fr(x^+_i)+\Fr(x^-_i)+y^+_i+y^-_i\,.
\end{align*}
By definition of the projections $\pi_j$, we can write
\begin{align*}
\pi_j(\mutn) = \mathcal{C}( \mutnIj, \lapj, \lamj, \gapj,\gamj ) =     \mutn_{|I} + \sum_{i=1}^{N^+\wedge j} \gap_i \delta_{\lap_i} + \sum_{i=1}^{N^-\wedge j} \gam_i \delta_{\lam_i}
\end{align*}
with a continuous $\mathcal{C}$, defined by
\begin{align*}
\mathcal{C}(\mu,x^+,x^-,y^+,y^-) = \mu + \sum_{i=1}^{j} y^+_i \delta_{x^+_i} + \sum_{i=1}^{j} y^-_i \delta_{x^-_i}
\end{align*}
Note that $\mathcal{C}$ is not a bijection: point masses in $I$ may come from $\mu$ or from the points $x^+,x^-$. However, for a given $\tilde{\mu} \in \mathcal{S}$, we may still easily calculate
\begin{align*}
\tilde{\mathcal{I}}_j(\tilde{\mu})= \inf \left \{ \mathcal{I}(\mu,x^+,x^-,y^+,y^-) |\, \mathcal{C}(\mu,x^+,x^-,y^+,y^-) = \tilde{\mu} \right\}
\end{align*}
as the minimum is attained by choosing $\mu=\tilde{\mu}_{|I(j)}$ and $\lambda_i^+=\ap$ and $\gamma_i^+=0$ if $i>N^+$ and $\lambda^-_i=\am$ and $\gamma_i^-=0$ if $i>N^-$. Therefore
\begin{align*}
\tilde{\mathcal{I}}_j(\tilde{\mu})=  \mathcal{K}( \mu_V \!\ | \!\ \tilde{\mu}) +\tilde{\mu}(I)- 1+ \sum_{i=1}^{N^+\wedge j} \left(\Fr(\lambda^+_i)+\gamma_i^+\right)+\sum_{i=1}^{N^-\wedge j}\left(\Fr(\lambda^-_i)+\gamma_i^-\right)
\end{align*}
which, by the contraction principle, is the rate function of $\pi_j(\mutn)$.
\QED

\subsection{Projective limit and normalization}
\label{susectPROJL}

From Theorem \ref{projLDP}, the projective method of the Dawson-G\"artner theorem, p. 162 in the book of \cite{demboz98}, yields the LDP for $\mutn$ under $\widetilde{\mathbb Q_n^V}$ with speed $\beta' n$ and rate function
\begin{align} \label{susectPROJLeq}
\tilde{\mathcal{I}}(\tilde\mu) = \sup_j \tilde{\mathcal{I}}_j(\tilde\mu) =  \mathcal{K}(\mu_V\ | \ \tilde\mu) +\tilde\mu(\R)- 1+\sum_{i=1}^{N^+} \Fr(\lambda^+_i)+\sum_{i=1}^{N^-}\Fr(\lambda^-_i) ,
\end{align}
defined for $\tilde\mu \in \Sr$. Recalling (\ref{normunnorm}), we want to come back to a normalized measure $\mu \in \Sr_1$. It would be natural to apply the mapping $\tilde\mu \mapsto \tfrac{\tilde\mu }{\tilde\mu(\R)}$ but unfortunately, this mapping is not continuous in our topology induced by \eqref{strangetop}. As a workaround, note that from the LDP for $\pi_j(\mutn)$, we also get the joint LDP of 
\begin{align*}
\big(\pi_j(\mutn) , \pi_j(\mutn)(\R) \big)
\end{align*}
as the mapping 
\begin{align*}
\pi_j(\tilde\mu) \longmapsto  \pi_j(\mutn)(\R) = \mutn(I) + \sum_{i=1}^{N^+\wedge j} \gamma_i^+ + \sum_{i=1}^{N^-\wedge j}\gamma_i^- 
\end{align*}
is continuous in our topology for any $j$. Thus, applying the projective method to $(\pi_j(\mutn) , \pi_j(\mutn)(\R) )$, we get the LDP for the pair $(\mutn,\mutn(\R))$ with rate function
\begin{align*}
\overline {\mathcal{I}}(\tilde\mu,\kappa) =  \begin{cases}  \tilde{\mathcal{I}}(\tilde\mu)\;\;\;& \mbox{if}\ \tilde\mu(\R)=\kappa\\ 
  \infty\;\; & \mbox{otherwise.}
\end{cases}
\end{align*} 
 Now we are able to recover the original spectral probability measures $\mun$ from the unnormalized measures $\mutn$, by applying to the pair $(\mutn,\mutn(\R))$ the (continuous) mapping $(\tilde\mu,\kappa) \mapsto \kappa^{-1}\tilde\mu$. The contraction principle yields then the LDP for $\mun$ under $\widetilde{\mathbb Q_n^V}$ (hence under $\mathbb Q_n^V$) with rate function
\begin{align*}
\mathcal{I}(\mu) = \inf_{\nu=\kappa \cdot \mu,\, \kappa>0} \tilde{\mathcal{I}}(\nu) = \inf_{\kappa>0} \tilde{\mathcal{I}}(\kappa \cdot \mu)
\end{align*}
By \eqref{susectPROJLeq}, we need to minimize over $\kappa$ the function
\begin{align*}
  - \int\log\left(\kappa\frac{d\mu}{d\mu_V}\right)d\mu_V  - 1+\kappa+ \sum_{i=1}^{N^+} \Fr(\lambda^+_i)+\sum_{i=1}^{N^-}\Fr(\lambda^-_i).
\end{align*}
The term $\kappa - 1 - \log \kappa = L^*(\kappa)$ attains its minimal value $0$ for $\kappa = 1$. We obtain therefore the following LDP.

\medskip

\begin{thm} \label{thmsusectPROJL}
The sequence of spectral measures $\mun$ under $\mathbb Q_n^V$, as a random element of $\Sr_1$ equipped with the topology induced by \eqref{strangetop}, satisfies the LDP with speed $\beta'n$ and rate function
\begin{align*}
\mathcal{I}_V(\mu) = \mathcal{K}( \mu_V \!\ | \!\ \mu) + \sum_{i=1}^{N^+} \Fr(\lambda^+_i)+\sum_{i=1}^{N^-}\Fr(\lambda^-_i).
\end{align*}
\end{thm}

\subsection{Topological considerations}
\label{topsection}

It remains to show that the LDP in Theorem \ref{thmsusectPROJL} holds in the weak topology and can be extended from $\Sr_1$ to $\mathcal{P}_1$. This step is a consequence of the following two lemmas. Their proofs are postponed to the appendix.

\begin{lem}\label{top1}
The weak topology on $\Sr_1$ is coarser than the topology induced by \eqref{strangetop}.
\end{lem}

\medskip

\begin{lem}\label{top2} The function $\mathcal{I}_V$, extended to $\mathcal{P}_1$ by setting $\mathcal{I}_V(\mu)=\infty$ if $\mu \notin \Sr_1$, is lower semicontinuous in the weak topology. 
\end{lem}

\medskip

By Lemma \ref{top1}, the LDP of Theorem \ref{thmsusectPROJL} holds with the weak topology on $\Sr_1$, and by Lemma \ref{top2}, the LDP can by extended to $\mathcal{P}_1$. This completes the proof of Theorem \ref{MAIN}.


\section{Appendix 1: Proof of Theorem \ref{LDPjextreme}}
\label{Appendix1}
\subsection{Comments on the assumptions}
\label{comm}
In \cite{benaych2012large}, the result is proved under their assumption 2.9, which consists of three requirements: confinement, technical condition and convergence of extreme eigenvalues. We may also refer to \cite{Benal} for a general result in the same vein. For the sake of completeness let us shortly discuss the three requirements.
\begin{enumerate}

\item Confinement

\begin{itemize}
\item[(B1)] 
$\displaystyle \qquad
\liminf_{|x| \rightarrow \infty} \frac{V(x)}{2\log |x|}> 1
 \, $
\end{itemize}

which differs from (A1) in the case $\beta < 1$. 
Actually, the proof of  the LDP  for the empirical spectral distribution $\mu_\u\sn$ of \cite{agz} does require  $\max (\beta^{-1}, 1)$ instead of $1$ since it is used to warrant the finiteness of 
$\int e^{-\beta' V(x)}dx$. Recently,   \cite{serf} proved the LDP for  $\mu_\u\sn$ under the assumption
\begin{itemize}
\item[(S1)] 
$\displaystyle \qquad
\lim_{|x| \rightarrow \infty}\!\ V(x) -  2 \log |x| = \infty
 \,. $
\end{itemize}
The proof is completely different, using the notion of $\Gamma$-convergence.
Besides, with another method (carrying everything on the unit circle by the Cayley transform),   the LDP for $\mu_\u\sn$ was proved in \cite{hardy2012note} under the weak confinement assumption
\begin{itemize}
\item[(H1)] 
$\displaystyle \qquad
\liminf_{|x| \rightarrow \infty}\!\ V(x) - 2 \max (1, \beta^{-1}) \log |x| > -\infty
 \,. $
\end{itemize}
Of course, when $\beta \geq 1$
\[ (A1) = (B1) \Rightarrow (S1) \ \text{and} \ (H1)\,,\]
but when $\beta < 1$ and $\liminf \frac{V(x)}{2\log |x|} \in (1 ,\beta^{-1})$, (S1) is satisfied but not  (H1). 

Since we use several times arguments taken from the proof of the LDP for $\muun$ in \cite{agz}, we did not weaken our hypothesis (A1)  into (B1) to avoid a complete rewriting, though we conjecture that this $\beta^{-1}$ is an artefact.

\item Technical condition.
\begin{itemize}
\item[(B2)]  For every $p \geq 1$, the limit 
$\displaystyle \qquad \lim_{n\to \infty} \frac{1}{n} \log \frac{Z^{n-p}_{\frac{n}{n-p}V}}{Z^n_{nV}}
 \,\ $  exists . 
\end{itemize}

\item Convergence of extreme eigenvalues.
\begin{itemize}
\item[(AGZ)]  Under $\Pnv$, the largest (resp. lowest) eigenvalue converges to $\alpha_\pm$ almost surely\footnote{see footnote 1}.
\end{itemize}
\end{enumerate}
In \cite{agz}  the LDP for $\lambda_{\min}$ is proved under (B2) with $p=1$. Later,   
in an erratum\footnote{available online http://www.wisdom.weizmann.ac.il/$\sim$zeitouni/cormat.pdf.}, the authors claimed that the proof of the LDP needs  one more assumption: either a slight modification of (AGZ) (replace $\Pnv$ by $\mathbb P^n_{nV/(n-1)}$) or  (A3) (control of large deviations). Later again, in \cite{Borot} a proof for $p=1$ is given under assumption (A3) alone, without (AGZ). It is worthwhile to mention the papers \cite{borot-multi} and \cite{alice-et-al} on connected topics.

To update the proof of \cite{benaych2012large} with the tools of \cite{Borot} adapted to the case $p > 1$ and for the sake of completeness, we give now the detailed scheme. We will use three lemmas whose proofs are postponed.

The first statement is a fact often mentioned (for instance  \cite{agz} pp. 83-84 or \cite{benaych2012large} p.744) but (as far as we know) never checked explicitly.  We set it as a lemma and for which we give a complete proof in Section \ref{531}  for convenience of the reader.

\begin{lem}
\label{LDPtilde}
Let  $V$ be  a potential satisfying the confinement condition (A1) and let $r$ be a fixed integer.  If  $\mathbb P^n_{V_n}$ is the probability measure associated to the potential
$V_n= \frac{n+r}{n} V$, 
then the law of $\muun$ under $\mathbb P^n_{V_n}$ satisfies the LDP with speed $\beta' n^2$ with good rate function

\begin{equation}  
\label{entropyV}
\mu \mapsto \mathcal E(\mu) - \inf_\nu\mathcal E(\nu)
\end{equation}
where $\mathcal E$  is defined in (\ref{ratemuu}).
\end{lem}

\begin{lem}
\label{teknik}
If the potential $V$ is finite and continuous on a compact set and infinite outside, we have, for every $p \geq 1$
\begin{equation}
\lim_{n\to \infty} \frac{1}{n} \log \frac{Z_V^n}{Z_{\frac{n}{n-p}V}^{n-p}} = -\inf_{x_1, \dots ,x_k} \sum_{k=1}^p  \mathcal J_V (x_k) = -p \inf_x \mathcal J_V(x)\,.
\end{equation}
\end{lem}

\begin{lem}\label{cvproba}
Under Assumption (A1) and (A3), the largest (resp. lowest) eigenvalue converges in probability to $\alpha^+$ (resp. $\alpha^-$).
\end{lem}

\subsection{Proof}
\subsubsection{Outline}
First notice that $\mathcal I_{\lambda^\pm}$ is a good rate function: it is lower semicontinuous as proved in \cite{Borot} A.1.p.478.  From the same reference, 
$\Fr_{V}^+$ and $\Fr_{V}^-$  have compact level sets, so that  
$\mathcal I_{\lambda^\pm}$ has compact level set by the union bound. The exponential tightness 
proved below  implies that the weak LDP for the extreme eigenvalues will be sufficient for the statement of Theorem \ref{LDPjextreme}. Throughout this proof, we will assume $\ell =j$, the generalization is straightforward. The weak LDP will follow from the upper bound
\begin{align}
\label{ub}
\limsup_{n\to \infty}\ (\beta'n)^{-1} \log \Pnv (\lambda^\pm(j) \in F^+ \times F^-) \leq - \inf_{(x^+,x^-)\in F^+ \times F^-} \mathcal{I}_{\lambda^\pm}(x^+, x^-)
\end{align}
for sets $ F^+= F^+_1 \times\dots \times F^+_j$ and idem for $F^-$, which generate the topology on $\R^{2j}$ and from the lower bound
\begin{align}
\label{lb}
\liminf_{n\to \infty}\ (\beta'n)^{-1} \log \Pnv(\lambda^\pm(j) \in G) \geq - \mathcal{I}_{\lambda^\pm}(x^+, x^-) ,
\end{align}
 for open sets $G$ containing $(x^+,x^-)$.

\subsubsection{Exponential tightness}
\label{expotight}

We define the compact set
\begin{align} \label{restrictV}
H_M= \{x \in [-M,M]: V(x) \leq M \}
\end{align}
and have to show that for every $j$ 
\begin{align}
\label{expoj}
\limsup_{M \rightarrow \infty}\limsup_{n \rightarrow \infty} n^{-1}\log \Pnv\big( (\lamj, \lapj) \notin H_M^{2j} \big) = -\infty.
\end{align}
But as
\begin{align*}
\Pnv\big( (\lamj, \lapj) \notin H_M^{2j} \big) \leq \Pnv\big( \lambda_1^+ >M  \big) + \Pnv\big( \lambda_1^- <-M  \big) ,
\end{align*}
exponential tightness reduces to the case $j=1$ and further (by symmetry),  we only have to prove 
\begin{equation}
\label{exptight}
\limsup_{M \rightarrow \infty}\limsup_{n \rightarrow \infty} n^{-1}\log \Pnv\big( \lambda_1^+ >M  \big)
 = -\infty\,. 
\end{equation}
By exchangeability, we have, 
\begin{eqnarray*}\Pnv(\lambda_1^+ \geq M) = 
n
\Pnv(\lambda_1  \geq M \ \hbox{and}\ \lambda_k \leq \lambda_1 \ \hbox{for}\  k = 2, \dots, n) \leq n \Pnv(\lambda_1  \geq M)
\end{eqnarray*}
Now
\begin{equation}
\label{bigintegral}
\Pnv(\lambda_1  \geq M) \leq \frac{Z^{n-1}_V}{Z^n_V}\int_{\{ x  \geq M\} }
 e^{-n\beta'V(x)}
\int
 \prod_{k=2}^n \left(|x-\lambda_k|^\beta e^{-\beta' V(\lambda_k)}\right) d\mathbb P^{n-1}_V(\lambda_{2}, \dots ,\lambda_n) dx
\end{equation}
Since $y \mapsto y^\beta$ is convex on $[0, \infty)$ for $\beta > 1$ and concave and subadditive for $\beta < 1$, we have
\[|x-\lambda|^\beta \leq (|x| + |\lambda|)^\beta \leq a(\beta)(|x|^\beta + |\lambda|^\beta) \]
with   $a(\beta) = 2^{(\beta -1)_+}$, and then
 \[|x-\lambda|^\beta e^{-\beta'V(\lambda)} \leq a(\beta)(|x|^\beta + |\lambda|^\beta) e^{-\beta' V(\lambda)}\,.\]
From the confinement assumption (A1) there exists $\eta > 0$ such that
\[\liminf_{|\lambda|\to \infty} \frac{V(\lambda)}{2\log |\lambda|} = 1 + 2\eta \]
and then there exists $M_0$ such that for $|\lambda| \geq M_0$ we have
\begin{equation}
\label{confin2}\frac{V(\lambda)}{2\log |\lambda|} \geq 1 + \eta\end{equation}
so that
\begin{eqnarray}
\label{comparison}
|\lambda|^\beta e^{-\beta'V(\lambda)} \leq |\lambda|^{-\beta\eta}\end{eqnarray}
(for $|\lambda| \geq M_0$) 
and 
then $|\lambda|^\beta e^{-\beta' V(\lambda)}$ is bounded by a constant $M_1$ uniformly in $\lambda$.
Now,  $V$ is bounded from below, say by $M_2$, and then
\[|x-\lambda|^\beta e^{-\beta'V(\lambda)} \leq a(\beta)(e^{-M_2} x^\beta + M_1)\]
and since $1 \leq x^\beta M_0^{-\beta}$ for $x\geq M_0$, we get
\[|x-\lambda|^\beta e^{-\beta'V(\lambda)} \leq b(\beta) x^\beta \leq b(\beta) e^{\frac{\beta'}{1+ \eta}V(x)}\]
with $b(\beta) = a(\beta)( e^{-M_2} + M_1 M_0^{-\beta})$
(the last inequality follows from (\ref{comparison})). 

Plugging this bound into (\ref{bigintegral}) we get
\[\Pnv(\lambda_1 \geq M) \leq \frac{Z_V^{n-1}}{Z_V^n}[b(\beta)]^{(n-1)}\int_M^\infty\exp \left(-\beta' \left(\frac{\eta n }{1+ \eta}\right)V(x)\right)\!\ dx\]
Since for $n$ large, $\beta'\eta n> (1+ \eta)/2$ we may write 
\[\int_M^\infty\exp \left(-\beta' \left(\frac{\eta n }{1+ \eta}\right)V(x)\right)\!\ dx \leq \exp  \left(- \left(\beta'\frac{\eta n }{1+ \eta}- \frac{1}{2}\right)C(M)\right)\times \int_M^\infty e^{-\frac{V(x)}{2}}  dx\]
where $C(M) = \inf \{ V(x) : |x| > M\}$. This last integral is finite in view of (\ref{confin2}).

We need the following lemma.
\begin{lem}
\label{Zratio}
\[\limsup_{n\to \infty} n^{-1} \log \frac{Z^{n-1}_V}{Z^n_V} := c_1 < \infty\,.\]
\end{lem}
Assuming the result of this lemma, we may write
\[\limsup_{n\to \infty} n^{-1}\log \Pnv(\lambda^+_1 \geq M) = \limsup_{n\to \infty} n^{-1} \log \Pnv (\lambda_1 > M)  \leq c_1 + \log b(\beta) - C(M)\frac{\beta'\eta}{1+\eta}\]
 which ends the proof of exponential tightness since $C(M) \rightarrow \infty$ when $M \rightarrow \infty$.

The proof of Lemma \ref {Zratio} is in \cite{Borot} p. 478 and makes use of the exponential tightness 
 of $(\mu_\u\sn)_n$ under assumption  (A1).


\subsubsection{Upper bound}
From the above paragraph, we may assume that all the $F^\pm_k$ are subsets of $H_M$. Furthermore, it suffices to consider $F_k^+$ subsets of $[\alpha^+,M]\cap H_M$ and $F_k^-$ subsets of $[-M,\alpha^-]\cap H_M$, as an extreme eigenvalue contained in $(\alpha^-+\varepsilon,\alpha^+-\varepsilon)$ implies that the distance of $\muun$ to $\mu_V$ in the weak topology is at least some $\delta>0$. From the LDP for $\muun$ with speed $n^2$ (Theorem \ref{LDPmuu}), 
 this event has negligible probability on our scale of speed $n$. \\
As before, we write 
$x^\pm$ for the vector $(x_1^+,\dots ,x_j^+,x_1^-,\dots ,x_j^-)$. After permutation, we may assume that the extreme eigenvalues are not among the eigenvalues $\lambda_0=(\lambda_1,\dots ,\lambda_{n-2j})$. Again to simplify notation and to omit several indicator functions, we will not assume throughout this proof that $x_1^+,\dots ,x_j^+$ or $x_1^-,\dots x_j^-$ are in the right order. The assertion for the ordered eigenvalues follows directly from the contraction principle.

We have the representation:
\begin{eqnarray}\label{representation}  \Pnv( \lambda^\pm(j) \in F^+\times F^-) = \frac{1}{Z^n_V} \frac{n!}{(j!)^2(n-2j)!} \int_{F^+ \times F^-} 
 \Upsilon_{n,j}(x^\pm)\ dx^\pm \end{eqnarray}
where
\begin{eqnarray}
\Upsilon_{n,j}(x^\pm)=  H(x^\pm)
 \Xi_{n,j}(x^\pm)  e^{-\beta'n\sum_1^{2j} V(x^\pm_k) }
\end{eqnarray}
with
\[H(x^\pm) =\prod_{1\leq r<s\leq 2j} |x^\pm_r - x^\pm_s|^\beta  \]
and, setting $\Delta(x^\pm) = (\max x_i^-, \min x_i^+)^{n-2j}$,
\begin{eqnarray}
\nonumber
\Xi_{n,j}(x^\pm) = \int_{\Delta(x^\pm) } \prod_{r=1}^{2j}   \prod_{s=1}^{n-2j} |x^\pm_r-\lambda_s|^\beta  \prod_{r=1}^{n-2j} e^{-n\beta'V(\lambda_r)}\prod_{1 \leq r <s \leq n-2j} |\lambda_r - \lambda_s|^\beta  d\lambda_0
\end{eqnarray} 
For $M$ large enough we may replace $V$ by $V_M$ 
defined by
\[V_M = \begin{cases}
V(x) \;\;\;&\mbox{if}\; x\in H_M,\\
\infty\;\;&\mbox{otherwise.}
\end{cases}\]
without any change (note that $V$ is necessarily bounded on $[\alpha^-,\alpha^+]$).
We have then
\begin{equation}\label{defxi}\Xi_{n,j}(x^\pm) = 
Z^{n-2j}_{\frac{n}{n-2j}V_M}
\int_{\Delta (x^\pm) }  \prod_{r=1}^{2j} \prod_{s=1}^{n-2j} |x^\pm_r-\lambda_s|^\beta   
 d \mathbb P^{n-2j}_{\frac{n}{n-2j}V_M}
(\lambda_0)
\end{equation} 
Set finally
\[Y_{n,n-2j}^M = \frac{Z^n_V}{Z^{n-2j}_{\frac{n}{n-2j}V_M}}\,.\]

We first find an upper bound for $\Upsilon_{n,j}(x^\pm)$. 
Let $B_\kappa =\{ \lambda_0 \in \R^{n-2j} : d(\mu_\u^{(n-2j)} , \mu_V) \leq \kappa\}$. On $\Delta(x^\pm)$ the integrand in \eqref{defxi} is bounded by $ e^{c_1n}$ for some $c_1=c_1(M) \geq 0$, so that
\begin{align}\label{proviso}
\left(Z^{n-2j}_{\frac{n}{n-2j}V_M}\right)^{-1}\Xi_{n,j}(x^\pm) & \leq \int_{\Delta(x^\pm) \cap B_\kappa }  \prod_{r=1}^{2j} \prod_{s=1}^{n-2j} |x^\pm_r-\lambda_s|^\beta   
 d\mathbb P^{n-2j}_{\frac{n}{n-2j}V_M}(\lambda_0) \\
& \notag \quad + e^{c_1n}  \mathbb P^{n-2j}_{\frac{n}{n-2j}V_M}(B_\kappa^c)\,.
\end{align}
One has to use Lemma \ref{LDPtilde}.
Since the rate function of the LDP for $\mu_\u^{(n-2j)}$ has a unique minimizer, Lemma \ref{LDPtilde} yields for the second term in the bound (\ref{proviso})
\begin{align*}
e^{c_1n}\!\  \mathbb P^{n-2j}_{\frac{n}{n-2j}V_M}(B_\kappa^c) \leq c_2e^{-c_3n^2}
\end{align*}
for some positive constants $c_2,c_3$. The right hand side of the first line in (\ref{proviso}) is bounded by
\begin{align*}
\exp \left\{ \beta (n-2j) \sup_{\mu:\ d(\mu,\mu_V)\leq \kappa}  \sum_{r=1}^{2j} \int \log |x^\pm_r-\eta| d\mu(\eta) \right\}
\end{align*} 
and then
\begin{align}
\label{ups}
&\left(Z_{\frac{n}{n-2j}V_M}^{n-2j}\right)^{-1}\Upsilon_{n,j}(x^\pm)\\
\nonumber
&
\leq H(x^\pm) \exp\left\{\beta'n\left( - \sum_{r=1}^{2j} V(x^\pm_r) +2 \sup_{\mu:\ d(\mu,\mu_V)\leq \kappa}\int \log| x^\pm_r- \eta| d\mu(\eta)\right)\right\}
+ \ e^{-c_4n^2}\,.
\end{align}
Recall the expression (\ref{poteff}) of the effective potential
\[\mathcal J_V(x) = V(x) - 2 \int \log|x-\xi| d\mu_V(\xi)\]
and use the bound
\begin{eqnarray}
\label{usethebound}
\limsup_{\kappa\downarrow 0} \sup_{\xi \in F} \sup_{d(\mu, \mu_V) \leq\kappa} \left(2\int \log|\xi - \eta|d\mu(\eta) - V(x) 
\right)
 \leq - \inf_{\xi \in F} \mathcal J_V(\xi)\end{eqnarray}
(see  \cite{Borot} p. 480) we get, for any $\eta > 0$ and $n$ large enough
\begin{eqnarray}
\label{480}\left(Z_{\frac{n}{n-2j}V_M}^{n-2j}\right)^{-1}\sup_{x^\pm \in F^\pm} \Upsilon_{n,j}(x^\pm )\leq \exp \beta' n\left(\eta- \inf_{x^\pm \in F^\pm} \sum_{r=1}^{2j} 
  \mathcal J_V(x^\pm_r) \right)\,,\end{eqnarray}
and then, owing to (\ref{representation}) and \eqref{defxi}, we get for any $\eta > 0$
\begin{eqnarray}\nonumber
\limsup_{n\to \infty} \  (\beta'n)^{-1}\log \Pnv (\lambda^\pm (j) \in F^\pm ) \leq \beta'\eta - \beta' \inf_{x^\pm \in F^\pm} \sum_{r=1}^{2j}
  \mathcal J_V(x^\pm_r) - \liminf_{n\to \infty}  (\beta'n)^{-1} \log Y_{n,n-2j}^M\,.
\end{eqnarray}
so that, since $\eta$ is arbitrary
\begin{align}
\limsup_{n\to \infty} \  (\beta'n)^{-1} \log \Pnv (\lambda^\pm (j) \in F^\pm ) \leq  -  \inf_{x^\pm \in F^\pm} \sum_{r=1}^{2j}
  \mathcal J_V(x^\pm_r) - \liminf_{n\to \infty}  (\beta'n)^{-1} \log Y_{n,n-2j}^M\,. \label{prov}
\end{align}
It remains to find a lower bound for $Y_{n,n-2j}^M$. 
We start from
\begin{align*}
Y_{n,n-2j}^{M}&= \frac{Z^n_V}{Z^{n-2j}_{\frac{n}{n-2j}V_M}}\geq  \frac{Z^{n}_{V_M}}{Z^{n-2j}_{\frac{n}{n-2j}V_M}}
\end{align*}
and we use the result of the Lemma \ref{teknik}, noticing that that for $M$ large enough the equilibrium measure is still $\mu_V$ and also that
$\inf \mathcal J_V(x) = \inf \mathcal J_{V_M}(x)$ for $M$ large enough.

Coming back to (\ref{prov}) yields to the expected upperbound (\ref{ub}).\\

\subsubsection{Lowerbound for large deviations} 
We start from an open ball $B= B(\xi^\pm, \varepsilon)$ centered at $\xi^\pm \in \R^{2j}$ with radius $\varepsilon$ in the sup-norm. Without loss of generality we may assume that it is included in $(\alpha^+,M)^j\times (-M,\alpha^-)^j$ as well as in $\{x\in \mathbb{R}^{2j} |\, V(x_i)\leq M \text{ for all } i\}$. We have again
\[Z^n_V\!\ \Pnv(\lambda^\pm (j) \in B)= \int_{B}  \Upsilon_{n,j}(x^\pm) dx^\pm \]

Let us consider the probability measure $\chi_j^M$ on $\R^n$, defined by
\begin{align*}
d \chi_j^M(x^\pm, \lambda) := (\kappa_{n,j}^M)^{-1} \mathbbm{1}_{ B}(x^\pm) \mathbbm{1}_{\Delta(x^\pm)}(\lambda_0) dx^\pm  
d\mathbb P^{n-2j}_{\frac{n}{n-2j}V_M}(\lambda_0) \end{align*}
where $\kappa_{n,j}^M$ is the normalizing constant. 
We have
\begin{eqnarray*}
\int_{B}  \Upsilon_{n,j}(x^\pm) dx^\pm = Z^{n-2j}_{\frac{n}{n-2j}V_M} \kappa_{n,j} I_{n,j}^M
\end{eqnarray*}
where
\begin{eqnarray*}
I_{n,j}^M := \int 
H(x^\pm)  e^{-\beta'n\sum_{k=1}^{2j} V_M(x^\pm_k)}
 \left(\prod_{k=1}^{n-2j}|x^\pm_r - \lambda_k|^\beta\right) d \chi_j^M(x^\pm, \lambda) 
\end{eqnarray*}
Jensen's inequality gives,
\begin{eqnarray*}
\frac{1}{\beta'}\log I_{n,j}^M
\geq n I_n^{(1)} + 2I_n^{(2)} + 2 (n-2j)I_n^{(3)} ,
\end{eqnarray*}
where
\begin{align*}
&I_n^{(1)} = -\int \sum_{k=1}^{2j} V_M(x^\pm_k) d \chi_j^M(x^\pm, \lambda),\\
&I_n^{(2)} = \int \sum_{1\leq r<s \leq 2j} \log |x^\pm_r-x^\pm_s| d \chi_j^M(x^\pm, \lambda),\\
&I_n^{(3)} = \frac{1}{n-2j}\int \sum_{r=1}^{n-2j}\sum_{k=1}^j \log |x^\pm_k - \lambda_r| d \chi_j^M(x^\pm, \lambda)\,.
\end{align*}

Lemma \ref{cvproba} implies that for $x^\pm \in B$ (recall that points in $B$ are bounded away from the support of $\mu_V$), 
\begin{eqnarray}
\label{lim1}
\mathbb P^{n-2j}_{\frac{n}{n-2j}V_M}(\lambda_0 \in \Delta(x^\pm) ) \xrightarrow[n \rightarrow \infty ]{}  1
\end{eqnarray}
and then 
\begin{equation}
\label{kappainfty}\kappa_{n,j}^M \xrightarrow[n \rightarrow \infty ]{} \int_B dx^\pm = (2\varepsilon)^{2j} \end{equation}
We have that
$\sum_k V(x^\pm_k)$ is uniformly bounded  on $B$. So that,  from (\ref{lim1}) and (\ref{kappainfty})
\begin{eqnarray}
\lim_{n \to \infty} \frac{1}{\beta'n}  I_n^{(1)} =
(2\varepsilon)^{-2j}\int_B \sum_{k=1}^{2j} V(x^\pm_k)\!\ dx^\pm\,.
\end{eqnarray}
For the second term we may assume without loss of generality that $B$ is tie-free.  Hence, we get a similar conclusion
\begin{eqnarray}
\lim_{n \to \infty} I_n^{(2)}
=  (2\varepsilon)^{-2j} \int_B \sum_{1\leq r<s \leq 2j} \log |x^\pm_r-x^\pm_s|\!\ dx^\pm\,.
\end{eqnarray}
Now for the third term, we first bound by below $ \mathbbm{1}_{z >0} \log z$  by $\log_M (z) :=(\log z) \mathbbm{1}_{0 < z\leq M}$ and set
\[ \ell_{M,B}(t) :=\int_{B}\sum_{k=1}^{2j}\log_M |x^\pm_k- t|\!\ dx^\pm \,,\]
so that
\[I_n^{(3)} \geq 
(\kappa_{n,j}^M)^{-1}\int  \int \ell_{M,B}(z) d\mu_\u^{(n-2j)}(z) 
d\mathbb P^{n-2j}_{\frac{n}{n-2j}V_M}
(\lambda_0)\,.\]
Since $\ell_{M,B}$ is  continuous and bounded as the convolution of a $L^1$ and a $L^\infty$ function, the above bound
converges to $(2\varepsilon)^{-2j} \int \ell_{M,B}(z) d\mu_V(z)$ as $n \rightarrow \infty$. Let us notice that since the support of $\mu_V$ is compact, and $B$ is fixed, we can choose $M$ large enough so that $x-M\leq t\leq x +M$ for every $x \in B$ and $t$ in the support of $\mu_V$. We get $\int \ell_{M,B}(z) d\mu_V(z) = \int \ell_B(z)  d\mu_V(z)$ where
\[\ell_B(t) =\int_{B}\sum_{k=1}^{2j}\log |x^\pm_k- t|  dx^\pm\,.\]
At this stage, we have :
\begin{eqnarray*}\liminf_{n\to \infty} \  \frac{1}{\beta'n} \log I_{n,j}^M \geq - (2\varepsilon)^{-2j} \int_B \sum_{k=1}^{2j} V(x^\pm_k) dx^\pm 
+ 2 (2\varepsilon)^{-2j} \int \ell_{B}(z) d\mu_V(z)\,.
\end{eqnarray*}
Going back to $\Pnv(\lambda^\pm (j) \in B)$ we get
\begin{align}
\notag
\liminf_{n\to \infty}  \frac{1}{\beta'n} \log \Pnv(\lambda^\pm (j) \in B) &\geq - \limsup_{n\to \infty}  \frac{1}{\beta'n} \log \frac{Z^n_V}{Z^{n-2j}_{\frac{n}{n-2j}V_M}}-  (2\varepsilon)^{-2j} \int_B \sum_{k=1}^{2j} V(x^\pm_k) dx^\pm\\
\label{goingback}
& \quad  + 2  (2\varepsilon)^{-2j}\int \ell_{B}(z) d\mu_V(z)\,.
\end{align}
Splitting into two parts the integral defining $Z_V^n$ we have
\[Z_V^n = 
Z_{V_M}^n
+ Z_V^n \!\ \Pnv ( \lambda_1^\pm\notin H_M)\]
and from the exponential tightness,
\[\frac{Z_V^n}{Z_{nV_M/n-p}^{n-p}} \leq  \frac{1}{1 - e^{-nC(M)}}\frac{Z_{V_M}^n}{Z_{nV_M/n-p}^{n-p}} \]
Now we apply Lemma \ref{teknik}
\[\limsup_{n\to \infty}  \frac{1}{\beta'n} \log \frac{Z^n_V}{Z^{n-2j}_{\frac{n}{n-2j}V_M}}
 \leq - \inf_{x^\pm} \sum_{k=1}^{2j} 
  \mathcal J(x^\pm_k)\]
which, plugged into (\ref{goingback}), yields
\begin{align}
\notag \liminf_{n\to \infty} \frac{1}{\beta'n} \log \Pnv(\lambda^\pm (j) \in B) &\geq \inf_{x^\pm} \sum_{k=1}^{2j} 
  \mathcal J(x^\pm_k) -  (2\varepsilon)^{-2j}\int_B \sum_{k=1}^{2j} V(x^\pm_k) dx^\pm \\
& \quad + 2 (2\varepsilon)^{-2j} \int \ell_{B}(z) d\mu_V(z) \,.
\end{align}
Remembering that $B$ has volume $(2\varepsilon)^{2j}$, and letting $\varepsilon \rightarrow 0$ we obtain
\begin{align*}
\liminf_{n\to \infty} \frac{1}{\beta'n} \log \Pnv(\lambda^\pm (j) \in B) \geq \inf_{x^\pm} \sum_{k=1}^{2j}   \mathcal J(x^\pm_k)-  \sum_{k=1}^{2j} V(\xi^\pm_k) + 2 \int \sum_{k=1}^{2j}\log |\xi^\pm_k- t|  d\mu_V(t) .
\end{align*}
This yields the lower bound (\ref{lb}) and completes the proof of Theorem \ref{LDPjextreme}.
\QED

\subsection{Proofs of Lemmas of Section 5.1}

\subsubsection{Proof of Lemma \ref{LDPtilde}}
\label{531}

Notice that
\begin{equation}
\frac{d\mathbb P^n_{V_n}}{d\Pnv} = \frac{Z^n_V}{Z^n_{V_n}} \exp \left( - \beta' rn  \muun(V)\right)
\end{equation}
We need the following lemma.
\begin{lem}
\label{zoverz}
\[\lim_{n\to \infty} \frac{1}{n^2} \log \frac{Z^n_{V_n}}{Z^n_V} = 0\]
\end{lem}
Admitting the result of this lemma, we follow the steps of  \cite{agz} Section 2.6.  \\

\textbf{Exponential tightness:}\\ We have
\begin{align*}
\mathbb P^n_{V_n} (\muun(V) > t)  &= \frac{ Z^n_V}{Z^n_{V_n}} 
\int_{\{ \muun(V) > t\} } \exp \left(- \beta' rn  \muun(V)  \right) d\Pnv \\
&\leq   \frac{ Z^n_V}{Z^n_{V_n}} e^{-n\beta' rt} \Pnv \left (\muun(V) > t)\right)
\end{align*}
which yields
\[\limsup_{n\to \infty} \frac{1}{n^2}\log \mathbb P^n_{V_n} (\muun(V) > t)) \leq \limsup_{n\to \infty} \frac{1}{n^2}\log  \Pnv (\muun(V) > t))\]
and we may refer to the classical case.\\

\textbf{Large Deviation upper bound:}\\
Let $\mu$ be any probability measure on $\R$. We start from
\begin{equation}
\label{basicRN} 
\mathbb P^n_{V_n}(d(\muun, \mu) \leq \varepsilon) =  \frac{ Z^n_V}{Z^n_{V_n}}\int_{\{ d(\muun, \mu) \leq \varepsilon \} } \exp \left(- \beta' rn  \muun(V)  \right) d\Pnv \end{equation}
Since $V$ is bounded below by $V_{\min}$, we get the upper bound
\begin{align*}
\mathbb P^n_{V_n}(d(\muun, \mu) \leq \varepsilon)\leq   \frac{ Z^n_V}{Z^n_{V_n}}
e^{-\beta' rn V_{\min}} \Pnv (d(\muun, \mu) \leq \varepsilon)
\end{align*}
so that
\[\limsup_{n\to \infty} \frac{1}{n^2}\log \mathbb P^n_{V_n}  (d(\muun, \mu) \leq \varepsilon)  \leq \limsup_{n\to \infty} \frac{1}{n^2}\log  \Pnv  (d(\muun, \mu) \leq \varepsilon) \]
and we may refer to the classical case.\\

\textbf{Large deviations lower bound:}\\
We start again from (\ref{basicRN}) and get, for every $t > 0$, the bound
\[\mathbb P^n_{V_n}(d(\muun, \mu) \leq \varepsilon) \geq \frac{ Z^n_V}{Z^n_{V_n}} 
e^{-\beta' rnt}\Pnv (d(\muun, \mu) \leq \varepsilon ,\, \muun(V) \leq t)\]
Now,
\begin{align*}
\Pnv (d(\muun, \mu) \leq \varepsilon ,\, \muun(V) \leq t) 
&\geq  \Pnv (d(\muun, \mu) \leq \varepsilon) - \Pnv ( \muun(V) > t)
\end{align*}
From the previous consideration of exponential tightness, it is possible to choose $t$ large enough  so that
\[\liminf_{n\to \infty} \frac{1}{n^2}\log \Pnv (d(\muun, \mu) \leq \varepsilon ,\, \muun(V) \leq t)  \geq \liminf_{n\to \infty}  \frac{1}{n^2}\log \Pnv (d(\muun, \mu) \leq \varepsilon)\]
and  refer to the classical case. This ends the proof of Lemma \ref{LDPtilde}.\\

{\bf Proof of Lemma \ref{zoverz}:}

We have 
\[\frac{Z^n_{V_n}}{ Z^n_V}= \int \exp \left( -\beta' rn \muun(V)\right)d\Pnv \]
On the one hand, we observe that since $V$ is bounded from below by a constant $V_{\min}$, we have
\[\frac{Z^n_{V_n}}{ Z^n_V} \leq \exp \left(-\beta' rn V_{\min}\right)\]
so that
\begin{equation}
\limsup_{n\to \infty} \frac{1}{n^2}\log \frac{Z^n_{V_n}}{ Z^n_V}
 \leq 0
\,.\end{equation}
On the other hand, for every $t > 0$ 
\[
\frac{Z^n_{V_n}}{ Z^n_V}
\geq e^{-\beta' rnt} \Pnv\left (\muun(V) \leq t\right) \]
From (2.6.21) in \cite{agz}, we know that  $\lim_{n\to \infty}  \Pnv\left (\muun(V) \leq t\right) = 1$  for $t$ large enough. We easily deduce 
\begin{equation}
\liminf_{n\to \infty}  \frac{1}{n^2}\log \frac{Z^n_{V_n}}{ Z^n_V}
 \geq 0
\,,
\end{equation}
which ends the proof of Lemma \ref{zoverz}.

\subsubsection{Proof of Lemma \ref{teknik}}

Let 
\[Y_{n,n-p} :=  \frac{Z_V^n}{Z_{\frac{n}{n-p}V}^{n-p}}\,.\]

\textbf{Lower bound:}\\
\begin{align*}
 \frac{Z_V^n}{Z_{\frac{n}{n-p}V}^{n-p}} &= \int \int \prod_{r=1}^{p} \left(  e^{-\beta' nV(x_r)}\prod_{k=1}^{n-p}|x_r - \lambda_k|^\beta\right) H(x) \!\ dx \!\ d \mathbb P^{n-p}_{\frac{n}{n-p}V}(\lambda_0)
\\
&\geq \int_{B(\xi, \varepsilon)} \int \prod_{r=1}^{p} \left(  e^{-\beta' nV(x_r)}\prod_{k=1}^{n-p}|x_r - \lambda_k|^\beta\right) H(x)\!\  dx\ \! 
 d \mathbb P^{n-p}_{\frac{n}{n-p}V}(\lambda_0)
\end{align*}
where $\xi$ is a point in $\mathbb R^p$ with tie-free entries and $B(\xi, \varepsilon)$ a ball of radius $\varepsilon$ in the sup-norm. We may assume that $B(\xi, \varepsilon)\subset \{x\in \mathbb{R}^{2j} |\, V(x_i)<\infty \text{ for all } i\}$. 
On $B(\xi, \varepsilon)$ the potential $V$ is uniformly continuous and we may replace each $x_k$ by $\xi_k$ and get for $\varepsilon$ small enough
\begin{align*}
 H(x) \prod_{r=1}^{p}   e^{-\beta' nV(x_r)} \geq \prod_{r=1}^{p}   e^{-\beta' nV(\xi_r)-\beta' n \delta_\varepsilon}
\end{align*}
for some $\delta_\varepsilon>0$. Set
 \[\ell_\varepsilon(y, \lambda) = (2\varepsilon)^{-1} \int_{y-\varepsilon}^{y+ \varepsilon} \log|t-\lambda| dt\]
Applying Jensen's inequality (the exponential is convex and we integrate over 
$(2\varepsilon)^{p} \prod_{r=1}^{p} \mathbbm{1}_{|x_r -\xi_r|\leq \varepsilon} dx_k$) we get 
\begin{align} \label{lowerbound1b}
Y_{n,n-p} 
&\geq  (2\varepsilon)^{p}  \left(\prod_{r=1}^{p}  e^{-\beta' nV(\xi_r) -\beta'n\delta_\varepsilon} \right)
\int \exp\left( 2\beta'(n-p)\sum_{r=1}^{p} \int \ell_\varepsilon(\xi_r, z)d
\mu_\u^{(n-p)}
 (z) \right) d\mathbb P^{n-p}_{\frac{n}{n-p}V} (\lambda_0) .
\end{align}
The function $\lambda \mapsto \ell_\varepsilon(y, \lambda)$ is  continuous and bounded (all variables live in a compact set).
 We can bound the integral in \eqref{lowerbound1b} from below by
\begin{align*}
& \quad \int_{\{ \lambda_0 : d(\mu_\u^{(n-2j)} , \mu_V) \leq \kappa\}} \exp\left( 2\beta'(n-p)\sum_{r=1}^{p} \int \ell_\varepsilon(\xi_r, z)d
\mu_\u^{(n-p)} (z)\right)  
d\mathbb P^{n-p}_{\frac{n}{n-p}V} (\lambda_0)\\
 &\geq 
\mathbb P^{n-p}_{\frac{n}{n-p}V}\left( d(\mu_\u^{(n-p)}, \mu_V)  \leq \varepsilon\right) \exp \left( 2 \beta'(n-p) \sum_{r=1}^{p}\int  \ell_\varepsilon(\xi_r, z) d\mu_V(z) -\beta' \delta'_\varepsilon n \right)
\end{align*}
with $\delta'_\varepsilon>0$, hence
\begin{align*}
Y_{n,n-p} \geq  (2\varepsilon)^{p} & \left( \prod_{r=1}^{p}  e^{-\beta' nV(\xi_r) -\beta'n\delta_\varepsilon} \right)
\mathbb P^{n-p}_{\frac{n}{n-p}V}
\left( d(
\mu_\u^{(n-p)}, \mu_V)  \leq \varepsilon \right)\\
&\times \exp \left( 2 \beta'(n-p) \sum_{r=1}^{p}\int  \ell_\varepsilon(\xi_r, z) d\mu_V(z) \right)
\end{align*}
According to Lemma \ref{LDPtilde}, we know that $\mathbb P^{n-p}_{\frac{n}{n-p}V}
\left( d(
\mu_\u^{(n-p)}
, \mu_V)  \leq \varepsilon\right)\rightarrow 1$. Since the 
 logarithmic potential $t\mapsto \int \log |t-\lambda| d\mu_V(\lambda)$ is continuous,
 we may 
write 
\[ \int  \ell_\varepsilon(\xi_r, z) d\mu_V(z) \geq \int \log |\xi_r -z|d\mu_V(z) - \delta''_\varepsilon\,,\]
where $\delta''_\varepsilon$ may depend on $\xi$ but tends to zero with $\varepsilon$.
We have then, for $\xi$ fixed, for every $\varepsilon > 0$
\begin{align*}
\liminf_{n\to \infty} \ (\beta' n)^{-1}\log Y_{n,n-p} \geq  
-\sum_{r=1}^{p} \mathcal J(\xi_r)
 -2p(\delta'_\varepsilon +\delta''_\varepsilon)
\end{align*}
Since it is true for every $\varepsilon$, after optimizing in $\xi$ we may conclude
\begin{eqnarray}
\liminf_{n \to \infty}\ (\beta'n)^{-1}\log Y_{n,n-p} \geq  - \inf_{\xi} \sum_{r=1}^{p} \mathcal J(\xi_r)\,.
\end{eqnarray}

\textbf{Upper bound:}\\
We have
\[Y_{n,n-p} = \int  H(x)\prod_{r=1}^p e^{-\beta' n \sum_1^p V(x_r)} \left(\int \prod_{k=1}^{n-p} |x_r - \lambda_k|^\beta \!\ d\mathbb P_{\frac{n}{n-p}V}^{n-p}(\lambda_0)\right)  dx\]
Recall the definition $B_\kappa =\{ \lambda_0 \in \R^{n-2j} : d(\mu_\u^{(n-2j)} , \mu_V) \leq \kappa\}$.  
Since all variables live on a compact set, $ \prod_{k} |x_r - \lambda_k|^\beta$  is bounded by $ e^{c_1n}$ for some $c_1 > 0$, and then
\begin{align}\label{provisob}
\int \prod_{k=1}^{n-p} |x_r - \lambda_k|^\beta \!\ d\mathbb P_{\frac{n}{n-p}V}^{n-p}(\lambda_0) \leq \int_{B_\kappa }   \prod_{k=1}^{n-p} |x_r - \lambda_k|^\beta  
 d\mathbb P^{n-p}_{\frac{n}{n-p}V}(\lambda_0)  + e^{c_1n}  \mathbb P^{n-p}_{\frac{n}{n-p}V}(B_\kappa^c)\,.
\end{align}
Since the rate function of the LDP  has a unique minimizer, Proposition \ref{LDPtilde} yields 
\begin{align*}
e^{c_1n}\!\  \mathbb P^{n-p}_{\frac{n}{n-p}V}(B_\kappa^c) \leq c_2e^{-c_3n^2}
\end{align*}
for some positive constants $c_2,c_3$. The integral on the right hand side of (\ref{provisob}) is bounded by
\begin{align*}
\exp \left\{ \beta (n-p) \sup_{\mu:\ d(\mu,\mu_V)\leq \kappa}  \sum_{r=1}^{p} \int \log |x_r-\eta| d\mu(\eta) \right\}
\end{align*} 
and then, since we integrate over a compact set,
\begin{align}
\label{ups1}
Y_{n,n-p}
&\leq \int  H(x) \exp\left\{\beta'n\left( - \sum_{r=1}^{p} V(x_r) +2 \sup_{\mu:\ d(\mu,\mu_V)\leq \kappa}\int \log| x_r- \eta| d\mu(\eta)\right)\right\} dx \\
\nonumber&\quad + \ c_4e^{-c_5n^2}\,.
\end{align}
If we 
 use again the bound (\ref{usethebound})
 we get, for any $\eta > 0$ and $n$ large enough
\begin{eqnarray}
\label{480b}
Y_{n,n-p}\leq \exp \beta' n\left(\eta- \inf_{x} \sum_{r=1}^{p} 
  \mathcal J_V(x_r) \right)\,,\end{eqnarray}
and then
\[\limsup_{n\to \infty} \frac{1}{\beta'n} \log Y_{n,n-p} \leq- \inf_{x} \sum_{r=1}^{p} 
  \mathcal J_V(x_r)\,.\]

\subsubsection{Proof of Lemma \ref{cvproba}}

From the LDP for extreme value of \cite{Borot},  the rate function is $\mathcal J_V - \inf_x \mathcal J_V(x)$. So, if this rate function which vanishes on the support of $\mu_V$ does not vanish outside, that means that the probability that $\lambda_1^+$ is greater than $\alpha + \varepsilon$ is exponentially small, and similarly for $\lambda_1^-$.

\section{Appendix 2 : Proof of Lemma \ref{top1} and  Lemma \ref{top2}}

\subsection{Proof of  Lemma \ref{top1}}
Let $\mu_n\rightarrow \mu $ 
in $\Sr_1$ 
equipped with 
 the topology induced by \eqref{strangetop}. Let $f$ be continuous and bounded and $\varepsilon>0$. Since $\mu$ is normalized, we 
may choose $N$ so large that 
\begin{align*}
\mu(I) + \sum_{i=1}^{N\wedge N^+} \gamma_i^+ + \sum_{i=1}^{N\wedge N^-} \gamma_i^- >1-\varepsilon . 
\end{align*}
Note that $N$ may be 0. Given this $N$, choose $n_0$ so large such that for all $n\geq n_0$
\begin{align*}
d_n := \left| \int g d\mu_{n|I} - \int g d\mu_{|I}\right| + \sum_{i=1}^{N\wedge N^+}| \gamma_{i,n}^+ g(\lambda_{i,n}^+) - \gamma_{i}^+ g(\lambda_{i}^+) | + \sum_{i=1}^{N\wedge N^-}| \gamma_{i,n}^- g(\lambda_{i,n}^-) - \gamma_{i}^- g(\lambda_{i}^-) | < \varepsilon 
\end{align*}
for $g\in \{1,f\}$, which is possible thanks to our topology on $\Sr$. 
This implies in particular
\begin{align*}
\sum_{i={N\wedge N^+}+1}^{N^+}| \gamma_{i,n}^+ | + \sum_{i={N\wedge N^-}+1}^{N^-}| \gamma_{i,n}^-| \leq 2\varepsilon .
\end{align*}
Then we have
\begin{align*}
& \qquad \left|\int f d\mu_{n} - \int f d\mu \right| \\
&\leq d_n + \sum_{i={N\wedge N^+}+1}^{N^+}| \gamma_{i,n}^+ f(\lambda_{i,n}^+)| + \sum_{i={N\wedge N^-}+1}^{N^-}| \gamma_{i,n}^- f(\lambda_{i,n}^-)| 
+ \sum_{i={N\wedge N^+}+1}^{N^+}| \gamma_{i}^+ f(\lambda_{i}^+)| +  \sum_{i={N\wedge N^-}+1}^{N^-}| \gamma_{i}^- f(\lambda_{i}^-)|  \\
&\leq d_n + 2\varepsilon ||f||_\infty + \varepsilon ||f||_\infty \leq \varepsilon + 3\varepsilon ||f||_\infty
\end{align*}
for all $n\geq n_0$.
\QED

\subsection{Proof of Lemma \ref{top2}}

We need to show that for measures $\mu_n \in \Sr_1$ with $\mu_n \to \mu \in \mathcal{P}_1\setminus \Sr_1$ weakly, we have $\mathcal{I}_V(\mu) \to \infty$. If $\mu \notin \Sr_1$, then either $\mu$ has a nondiscrete part outside of $[\alpha_-,\alpha_+]$ or an infinite number of atoms outside of $[\alpha_--\varepsilon , \alpha_++\varepsilon]$ for some $\varepsilon>0$. 

Now, the only way for $\mu_n$ to be arbitrarily close to such a measure $\mu$ in the weak topology is if there exists a number $\ell(n)$ of atoms $x_1,\dots ,x_{\ell(n)}$ that are not lying in $[\alpha_--\varepsilon , \alpha_++\varepsilon]$ for some $\varepsilon>0$ and $\ell(n) \to \infty$. This implies $\Fr(x_i)>\delta$ for some positive $\delta$ for all $i\leq \ell (n)$ and then $\mathcal{I}_V(\mu) \to \infty$.
\QED
\section*{Acknowledgments}
Many thanks are due both to Stanislas  Kupin and G\'erard Ben Arous for  helpful discussions at the time of the early genesis of this work. The spectral point of view 
 of Stanislas on sum rules has strongly enlighten the darkness of our research paths on this subject. During the conference for the birthday of one author of this work, G\'erard gave us the keys for the study for the large deviations of the spectral measure part  supported by the extreme eigenvalues. At the time of the first revision, we warmly thank Barry Simon for its insight remarks and suggestions, especially on the gems.
\bibliographystyle{apalike}
\bibliography{bibi}
\end{document}